\documentclass[a4paper]{article}
\usepackage{amssymb}
\usepackage{amsthm}
\usepackage{amsmath}
\usepackage{latexsym}
\usepackage[T1]{fontenc}
\usepackage[latin1]{inputenc}

\def \eop {\hbox{}\nobreak\hfill \vrule width 2.0mm height 1.8mm depth 0mm
\par \goodbreak \smallskip}

\newcommand{\integ}[2]{\displaystyle \int_{#1}^{#2}}
\newcommand{\dint}{\displaystyle\int}

\begin{document}
\numberwithin{equation}{section}
\newtheorem{definition}{Definition}[section]
\newtheorem{theorem}{Theorem}[section]
\newtheorem{proposition}{Proposition}[section]
\newtheorem{lemma}{Lemma}[section]
\newtheorem{remark}{Remark}[section]
\newtheorem{corollary}{Corollary}[section]
\newtheorem{example}{Example}[section]
\def \ce{\centering}
\def \bop {\noindent\textbf{Proof}}
\def \eop {\hbox{}\nobreak\hfill
\vrule width 2mm height 2mm depth 0mm
\par \goodbreak \smallskip}
\def \R{I\!\!R}
\def \N{I\!\!N}
\def \P{\mathbb{P}}
\def \E{I\!\!E}
\def \T{\mathbb{T}}
\def \H{\mathbb{H}}
\def \L{\mathbb{L}}
\def \Z{\mathbb{Z}}
\def \bf{\textbf}
\def \it{\textit}
\def \sc{\textsc}
\def \ni {\noindent}
\def \sni {\ss\ni}
\def \bni {\bigskip\ni}
\def \ss {\smallskip}
\def \F{\mathcal{F}}
\def \g{\mathcal{g}}
\def \eop {\hbox{}\nobreak\hfill
\vrule width 2mm height 2mm depth 0mm
\par \goodbreak \smallskip}
\title{Stochastic Quadratic BSDE With Two RCLL Obstacles}
\author{E. H. Essaky$^1$  \quad \quad M. Hassani$^1$ \quad \quad Y. Ouknine$^2$\\\\
$^1$ Universit\'{e} Cadi Ayyad\\ Facult\'{e} Poly-disciplinaire\\
D\'{e}partement de Math\'{e}matiques et d'Informatique\\ B.P. 4162,
Safi, Maroc.\\ e-mails : essaky@ucam.ac.ma
\hspace{.2cm}medhassani@ucam.ac.ma \hspace{.2cm} \\ \\
$^2$ Universit\'{e} Cadi Ayyad\\ Facult\'{e} des Sciences Semlalia\\
D\'{e}partement de Math\'{e}matiques\\ B.P. 2390,  Marrakech, Maroc.\\
e-mail : ouknine@ucam.ac.ma}
\date{}
\maketitle \maketitle \footnotetext[1]{This work is supported by
Hassan II Academy of Science and technology, Action Int\'egr\'ee
MA/10/224 and Marie Curie ITN n$^{\circ}$ 213841-2.}

\begin{abstract}
We study the problem of existence of solutions for generalized
backward stochastic differential equation with two reflecting
barriers (GRBSDE for short) under weaker assumptions on the data.
Roughly speaking we show the existence of a maximal solution for
GRBSDE when the terminal condition $\xi$ is ${\cal F}_T-$measurable,
the coefficient $f$ is continuous with general growth with respect
to the
 variable $y$ and stochastic quadratic growth with respect to the
 variable $z$  and the reflecting barriers $L$ and $U$ are
 just right continuous left limited. The result is proved without
 assuming any $P-$integrability conditions. 
\end{abstract}

\ni \textbf{Keys Words:} Backward stochastic differential equation;
stochastic quadratic growth; comparison theorem; exponential
transformation.
\medskip

\ni \textbf{AMS Classification}\textit{(1991)}\textbf{: }60H10,
60H20.
\bigskip

\section{Introduction}
Let $(\Omega, {\cal F}, ({\cal F}_t)_{t\leq T}, P)$ be a stochastic
basis on which is defined a Brownian motion $(B_t)_{t\leq T}$ such
that $({\cal F}_t)_{t\leq T}$ is the natural filtration of
$(B_t)_{t\leq T}$ and ${\cal F}_0$ contains all $P$-null sets of
$\cal F$.  Note that $({\cal F}_t)_{t\leq T}$ satisfies the usual
conditions, \it{i.e.} it is right continuous and complete.

The notion of BSDE with two reflecting barriers (RBSDE for short)
has been first introduced by Civitanic and Karatzsas \cite{CK}. A
solution for such an equation, associated with a coefficient $f$;
terminal value $\xi$ and two barriers $L$ and $U$, is a quadruple of
processes $(Y, Z, K^+, K^-)$ with values in $\R\times\R^d\times
\R_+\times \R_+$ satisfying:
\begin{equation}
\label{eq000} \left\{
\begin{array}{ll}
(i) &
 Y_{t}=\xi
+\integ{t}{T}f(s,Y_{s},Z_{s})ds +\integ{t}{T}dK_{s}^+
-\integ{t}{T}dK_{s}^- -\integ{t}{T}Z_{s}dB_{s}\,, t\leq T,
\\ (ii)&
   L_t \leq Y_{t}\leq U_{t},\quad \forall t\leq T,\\
(iii)& \integ{0}{T}( Y_{t}-L_{t}) dK_{t}^+= \integ{0}{T}(
U_{t}-Y_{t}) dK_{t}^-=0,\,\, P-\mbox{a.s.}, \\ (iv) &  K^+, K^-
\,\,\mbox{are continuous nondecreasing processes with }\,\, K_0^+
=K_0^- =0.
\end{array}
\right.
 \end{equation}
Here two continuous increasing processes $K^+$ and $K^-$ have been
added in order to force the solution $Y$ to remain in the region
enveloped by the lower reflecting obstacle $L$ and the upper
reflecting obstacle $U$. This is done by the
 cumulative actions of processes $K^+$ and $K^-$. In the case
of a uniformly Lipschitz coefficient $f$ and a square terminal
 condition $\xi$ the existence and uniqueness of a solution have been proved
 when the barriers $L$ and $U$ are either regular or satisfy Mokobodski's condition. This
last condition essentially postulates the existence of a
quasimartingale between the barriers $L$ and $U$. It has been shown
also in \cite{CK} that the solution coincides with the value of a
stochastic Dynkin game of optimal stopping. The link between
obstacle PDEs and RBSDEs has been given in Hamad\`ene and Hassani
\cite{HH}. 

The problem of existence of solutions for generalized BSDE with two
reflecting barriers under weaker assumptions on the input data has
been studied by Essaky and Hassani \cite{EH1} (see also \cite{EH2}
for the non-reflected case). The authors have proved the existence
of a maximal solution when the terminal condition $\xi$ is ${\cal
F}_T-$measurable, the coefficient $f$ is continuous with general
growth with respect to the
 variable $y$ and stochastic quadratic growth with respect to the
 variable $z$  and the reflecting barriers $L$ and $U$ are
 continuous. The result has been proved without assuming any $P-$integrability conditions.
Applications to the Dynkin game problem as well as to the American
game option have been also given.

In this paper we add a term of the form $\displaystyle{\sum_{t<s\leq
T}}h(s, Y_{s-}, Y_s)$ in the RBSDE (\ref{eq000}) where $h$ is a
process with values $\R$, which can be interpreted as a parameter of
jump reflection at barriers $L$ and $U$ (see Definition \ref{def1}).
A natural question is then arises : is there any solution for this
new RBSDE under the same assumptions as in \cite{EH1} but when the
barriers $L$ and $U$ are only right continuous left limited
(\it{rcll} for short) processes?. The present work gives a positive
answer to this question. The difficulty here lies in the fact that
since the barrier $L$ and $U$ are allowed to have jumps then the
process $Y$ is so and then the reflecting processes $K^+$ and $K^-$
are no longer continuous but just \it{rcll}. In this case, if $(Y,
Z, K^+, K^-)$ is a solution then its size of jumps is given by : 
$$
\begin{array}{ll}
 &  K^+_t -K_{t-}^+ =
\bigg(L_{t-}- [Y_t+h(t, Y_{t-}, Y_t )]\bigg)^+
\\ & K^-_t -K_{t-}^-  =
\bigg(U_{t-}- [Y_t+h(t, Y_{t-}, Y_t )]\bigg)^-
\\
& Y_{t-} = L_{t-} \vee [Y_t + h(t, Y_{t-}, Y_t)]\wedge U_{t-}.
 \end{array}
$$
By means of an exponential change, the proof of our main result
consists in establishing first a correspondence between our GRBSDE
and another GRBSDE whose coefficients are more tractable. We show
that the existence of solutions for our initial GRBSDE is equivalent
to the existence of solutions for the auxiliary GRBSDE. In order to
prove that this auxiliary GRBSDE admits a maximal solution the
following four cases are discussed :
\begin{enumerate}
\item  $f=h=0$.
\item $f$ is Lipschitz and $h=0$.
\item $f$ is Lipschitz and there exists a finite family of stopping times $S_0 =0\leq S_1\leq...\leq S_{p+1}=T$ such that
for every $x,y\in\R$ and $t\notin \{S_1,...,S_{p+1}\}$,\,\, $h(t,
\omega, x, y)=0.$
\item The general case.
\end{enumerate}
In the fourth case, since the integrability conditions on parameters
are weaker, we make use of approximations and truncations to
establish the existence result for the auxiliary GRBSDE. The final
step consists in justifying the passage to the limit and identifying
the limit as the solution of the auxiliary GRBSDE.

This paper is organized as follows. In the next section we lay out
the notation and the assumptions and state the main result. In
Section 3, by means of an exponential change, we show that the
existence of solutions for our initial GRBSDE is equivalent to the
existence of solutions for an auxiliary GRBSDE whose coefficients
are more tractable. Section 4 is devoted to the proof of our main
result. A comparison theorem for maximal solutions is proved in
Section 5. Finally, in the appendix we prove a comparison theorem
for solutions of GRBSDE which plays a crucial role in our proofs.

\section{Statements and main result for GRBSDE}
\subsection{Notations}
Let $(\Omega, {\cal F}, (\F_t)_{t\leq T}, P)$ be a stochastic basis
on which is defined a Brownian motion $(B_t)_{t\leq T}$ such that
$({\cal F}_t)_{t\leq T}$ is the natural filtration of $(B_t)_{t\leq
T}$ and ${\cal F}_0$ contains all $P$-null sets of $\cal F$.  Note
that $({\cal F}_t)_{t\leq T}$ satisfies the usual conditions,
\it{i.e.} it is right continuous and complete. For simplicity, we
omit sometimes dependence on $\omega$ of some processes or random
functions.
\medskip\\ Let us now introduce the following notations :
\medskip

$\bullet$ $\cal P$ the sigma algebra of ${\cal F}_t$-predictable
sets on $\Omega\times [0,T].$

$\bullet$ ${\cal D}$ is the set of ${\cal P}$-measurable and right
continuous with left limits (\it{rcll} for short) processes
$(Y_t)_{t\leq T}$ with values in $\R$.

$\bullet$ For a given process $Y\in {\cal D}$, we denote :
$Y_{t-}=\displaystyle{\lim_{s\nearrow t}Y_s}, t\leq T$
$(Y_{0^-}=Y_0)$ and $\Delta_s Y = Y_s - Y_{s-}$ the size of its jump
at $s$.

$\bullet$ ${\cal K} := \{ K\in D\quad :\quad K \quad\mbox{is
nondecreasing and }\,\, K_0  =0\}$.


$\bullet$ ${\cal K}^{c} := \{ K\in {\cal K}\quad :\quad \Delta_t K
=0,\,\, \forall t\in [0, T]\}$.

$\bullet$ ${\cal K}-{\cal K}$ the set of ${\cal P}$-measurable and
\it{rcll} processes $(V_t)_{t\leq T}$ such that there exist $V^+,
V^-\in{\cal K}$ satisfying : $V = V^+ -V^-$. In this case, for each
$\omega\in\Omega$, $dV_t(\omega)$ denotes the signed measure on
$([0, T], {\cal B}_{[0,T]})$ associated to $V_t(\omega)$ where
${\cal B}_{[0,T]}$ is the Borel sigma-algebra on $[0, T]$.

$\bullet$ For a given process $V\in {\cal K}-{\cal K}$, we define :
$\dint_a^b dV_s = V_b -V_a = \dint_{]a, b]}dV_s$ and $V_t^c = V_t
-\displaystyle{\sum_{0 <s\leq t}}\Delta_s V$.

$\bullet$ ${\cal L}^{2,d}$ the set of $\R^d$-valued and $\cal
P$-measurable processes $(Z_t)_{t\leq T}$ such that
$$\integ{0}{T}|Z_s|^2ds<\infty, P- a.s.$$   \\
%

\ni The following notations are also needed :\\

$\bullet$ For a stopping time $\nu$, \,\, $ [|\nu|] := \{(t,
\omega)\in [0, T]\times \Omega\,\, :\,\, \nu(\omega) = t\}.$

$\bullet$ For a set $B$, we denote by $B^c$ the complement of $B$
and $1_B$ denotes the indicator of $B$.

$\bullet$ For each $(a, b)\in\R^2$, $a\wedge b = \min(a, b)$\,\, and
\,\,$a\vee b = \max(a, b)$.

$\bullet$ For all $(a, b, c)\in\R^3$ such that $a\leq c$,\,\, $a\vee
b\wedge c = \min(\max(a,b), c) = \max(a, \min(c,b))$.



\subsection{Definitions}
Throughout the paper we introduce the following data :\\

\ni$\bullet$ $\xi$ is an ${\cal F}_T$-measurable one dimensional
random
variable.\\ 

\ni$\bullet$ $L:=\left\{ L_{t},\,0\leq t\leq T\right\}$ and
$U:=\left\{ U_{t},\,0\leq t\leq T\right\}$ are two barriers
 which belong to ${\cal D}$ such that $L_t\leq U_t$, $\forall t\in
 [0,T[$ and assume, without loss of generality, that $L_T= \xi= U_T$. \\

\ni $\bullet$ $f : [0,T] \times \Omega \times \R\times
\R^d\longrightarrow \R$ is a function  such that :
$$
\forall (y,z)\in \R\times \R^d,\,\,\, (t,\omega)\longmapsto f(t,
\omega, L_t(\omega)\vee y \wedge U_t(\omega), z)\,\,\mbox{ is }\,\,
{\cal P}-\mbox{measurable}.
$$


 \ni$\bullet$ $g :  [ 0,T]\times \Omega \times \R\longrightarrow \R$ is a
 function such that
$$
\forall y\in \R,\,\,\, (t,\omega)\longmapsto g(t, \omega,
L_t(\omega)\vee y \wedge U_t(\omega))\,\,\mbox{ is }\,\, {\cal
P}-\mbox{measurable}.
$$

\ni$\bullet$ $h : ]0,T] \times \Omega\times
\R\times\R\longrightarrow \R$ is a function such that
$$
\forall (t,x,y)\in ]0, T]\times\R\times \R,\,\,\, \omega\longmapsto
h(t, \omega,L_{t-}(\omega)\vee x \wedge U_{t-}(\omega)
,L_t(\omega)\vee y \wedge U_t(\omega))\,\,\mbox{ is }\,\,{\cal
F}_t-\mbox{measurable}.
$$


\ni$\bullet$ $A$ is a process in ${\cal K}^c$.\\

\medskip
\ni To give conditions under which solutions to a GRBSDE exist, we
should first give the following definitions.
\begin{definition}
Let $K^1$ and $K^2$ be two processes in ${\cal K}$. We say that :
\begin{enumerate}
\item $K^1$ and $K^2$ are singular if and only if there exists a set
$D\in {\cal P}$ such that
$$
\E\dint_0^T 1_D(s,\omega) dK^1_s(\omega) = \E\dint_0^T
1_{D^c}(s,\omega) dK^2_s(\omega) =0.
$$
 This is denoted by $dK^1 \perp dK^2$.
\item $dK^1 \leq dK^2$ if and only if for each set $B\in {\cal P}$
$$
\E\dint_0^T 1_B(s,\omega) dK^1_s(\omega)\leq \E\dint_0^T
1_{B}(s,\omega) dK^2_s(\omega), \quad \it{i.e.}\quad K_t^1
-K_s^1\leq K_t^2-K_s^2,\,\,\, \forall s\leq t\quad P-a.s.
$$
In this case $\dfrac{dK^1}{dK^2}$ denotes a ${\cal P}-$measurable
Radon-Nikodym density of $dK^1$ with respect to $dK^2$ which
satisfies
$$
0\leq \dfrac{dK^1}{dK^2}(s,\omega)\leq 1,\quad dK^2_s
(\omega)P(d\omega)-a.e. \,\,\mbox{on}\,\, [0,T]\times \Omega.
$$
\end{enumerate}
\end{definition}

Let us now introduce the definition of our GRBSDE with two \it{rcll}
obstacles $L$ and $U$.
\begin{definition}\label{def1}
\begin{enumerate}
 \item We say that $(Y,Z,K^+,K^-):=( Y_{t},Z_{t},K_{t}^+,K_{t}^-)_{t\leq T}$
is a solution of the GRBSDE, associated with the data $(\xi, f, g,
h, A, L, U)$, if the following hold :
\begin{equation}
\label{eq0} \left\{
\begin{array}{ll}
(i) & 
 Y_{t}=\xi
+\integ{t}{T}f(s,Y_{s},Z_{s})ds+\dint_t^Tg(s,
Y_s)dA_s +\sum_{t<s\leq T}h(s, Y_{s-}, Y_s)\\
&\qquad\quad+\integ{t}{T}dK_{s}^+ -\integ{t}{T}dK_{s}^-
-\integ{t}{T}Z_{s}dB_{s}\,, t\leq T,
\\ (ii)&
\forall t\in[0,T[,\,\, L_t \leq Y_{t}\leq U_{t},\\  (iii)&
\integ{0}{T}( Y_{t-}-L_{t-})
dK_{t}^+= \integ{0}{T}( U_{t-}-Y_{t-}) dK_{t}^-=0,\,\, \mbox{a.s.}, \\
(iv)& Y\in {\cal D}, \quad K^+, K^-\in {\cal K}, \quad Z\in {\cal
L}^{2,d},  \\ (v)& dK^+\perp  dK^-.
\end{array}
\right. \end{equation}
\item We say that the GRBSDE (\ref{eq0}) has a maximal (resp. minimal)
solution $(Y_t ,Z_t ,K^+_t , K_t^- )_{t\leq T}$ if for any other
solution $(Y_t^{'} ,Z_t^{'} ,K'^{+}_t , K'^{-}_t )_{t\leq T}$ of
(\ref{eq0}) we have for all $t \leq T$, $Y_t^{'}\leq Y_t$, $P$-a.s.
(resp. $Y_t^{'}\geq Y_t$, $P$-a.s.).
\end{enumerate}
\end{definition}
\begin{remark}In our definition we introduce
a process $h$ with values $\R$, which may be interpreted as a
parameter of jump reflection at barriers $L$ and $U$. Moreover, if
$(Y, Z, K^+, K^-)$ is a solution of GRBSDE (\ref{eq0}) then it
satisfies for all $t\in ]0, T]$ (see Lemma \ref{lem1})
$$
\begin{array}{ll}
 & \Delta_t K^+  =
\bigg(L_{t-}- [Y_t+h(t, Y_{t-}, Y_t )]\bigg)^+
\\ & \Delta_t K^-  =
\bigg(U_{t-}- [Y_t+h(t, Y_{t-}, Y_t )]\bigg)^-
\\
& Y_{t-} = L_{t-} \vee [Y_t + h(t, Y_{t-}, Y_t)]\wedge U_{t-}.
 \end{array}
$$
\end{remark}

\subsection{Assumptions and remarks}
We shall need the following assumptions on $f$, $g$, $h$, $L$ and
$U$:
\medskip

\ni $(\bf{A.1})$ There exist two processes $\eta \in L^0(\Omega,
L^1([0,T], ds, \R_+))$ and $C\in {\cal D}$ such that :
\begin{itemize}
 \item[\bf{(a)}]$\forall (y, z)\in
\R\times\R^d,\,\, |f(t, \omega, L_t(\omega)\vee y \wedge
U_t(\omega), z )| \leq \eta_t(\omega)+C_t(\omega)|z|^2,\quad
dtP(d\omega)-$a.e.
\item[\bf{(b)}] $dtP(d\omega)-$a.e.,\,\,the function $(y, z)\longmapsto f(t,
\omega, L_t(\omega)\vee y \wedge U_t(\omega), z)$ is continuous.
\end{itemize}\vspace{0.2cm}

\ni $(\bf{A.2})$ There exists $\beta \in L^0(\Omega, L^1([0,T],
A(dt), \R_+))$ such that :
\begin{itemize}
\item[\bf{(a)}]$
\forall y\in\R,\quad |g(t, \omega, L_t(\omega)\vee y \wedge
U_t(\omega))| \leq \beta_t(\omega),\quad A(dt)P(d\omega)-$a.e.
\item[\bf{(b)}] $A(dt)P(d\omega)-$a.e.,\,\,the function\,\,
$y\longmapsto g(t, \omega, L_t(\omega)\vee y \wedge U_t(\omega))$ is
continuous.
\end{itemize}\vspace{0.2cm}

\ni $(\bf{A.3})$
\begin{itemize} \item[\bf{(a)}] There
exists $l : ]0, T]\times \Omega \longrightarrow \R_+$ satisfying for
each $t\in ]0, T]$, $l_t$ is ${\cal F}_t$-measurable and
$P$-a.s.,$\displaystyle{\sum_{0<s\leq T}l_s} <+\infty$, such that:
$$
P-\mbox{a.s}., \forall (t, x, y)\in]0, T]\times \R\times\R,\quad
|h(t, \omega, L_{t-}(\omega)\vee x \wedge U_{t-}(\omega),
L_t(\omega)\vee y \wedge U_t(\omega))|\leq l_t(\omega).
$$
\item[\bf{(b)}] $P-$\mbox{a.s}., $\forall (t, y)\in ]0, T]\times \R,\,\,  x\longmapsto h(t, \omega,L_{t-}(\omega)\vee x \wedge
U_{t-}(\omega) ,L_t(\omega)\vee y \wedge U_t(\omega))$ is
continuous.

\item[\bf{(c)}] $P-$a.s., $\forall (t, x)\in]0, T]\times \R$,\,\, the function $y\mapsto y+h(t,
\omega, L_{t-}(\omega)\vee x \wedge U_{t-}(\omega) ,L_t(\omega)\vee
y \wedge U_t(\omega))$ is nondecreasing and continuous.
\end{itemize}

\noindent$(\bf{A.4})$ There exists a semimartingale $S_. = S_0 + V_.
+\dint_0^. \gamma_s dB_s$, with $S_0\in\R, V \in \cal K - \cal K$
and $\gamma\in{\cal L}^{2,d}$, such that $L_t \leq S_t\leq U_t, \,\,
\forall t\in [0,T].$\\

Before giving the main result of this paper, let us give the
following remarks on the assumptions.
\begin{remark}\label{rem1}
\begin{enumerate}
\item It should be pointed out that conditions $\bf{(A.1)(a)}$,
$\bf{(A.2)(a)}$ and $\bf{(A.3)(a)}$ hold if the functions $f$, $g$
and $h$ satisfy the following: $\forall (s,\omega),\,\,\forall y\in
[L_s(\omega), U_s(\omega)],\,\, \forall z\in \R^d,\,\,$ $\forall
x\in [L_{s-}(\omega), U_{s-}(\omega)]$
$$
\begin{array}{lll}
& |f(s, \omega, y, z )| \leq \widetilde{\eta}_s(\omega)\,\,\theta(s,
\omega, x, y) + \theta(s, \omega, x, y)|z|^2,
\\ &
|g(s, \omega, y)| \leq \widehat{\eta}_s(\omega)\,\,\theta(s,
\omega,x, y),
\\ &
|h(s, \omega,x, y)| \leq \widetilde{l}_s (\omega)\,\,\theta (s,
\omega,x, y),
\end{array}
$$
where:\\ $\bullet$  $\theta : [0,T]\times \Omega \times \R\times\R
\longrightarrow \R_+$ is a function such that for every $t\in [0,
T]$:
$$
P-a.s.,\,\, D_t(\omega) := \displaystyle\sup_{s\leq t, \alpha,
\delta\in [0,1]}\theta(s,\omega, \delta L_{s-}(\omega) +(1-\delta)
U_{s-}(\omega), \alpha L_s(\omega) +(1-\alpha) U_s(\omega)) <
+\infty,
$$
and ${\cal F}_t-$adapted.

$\bullet$ $\widetilde{\eta} \in L^0(\Omega, L^1([0,T], ds,
\R_+))$ and $\widehat{\eta} \in L^0(\Omega, L^1([0,T], dA_s, \R_+))$,\\
$\bullet$ $\widetilde{l}_t$ is an ${\cal F}_t$-measurable function
satisfying
$\displaystyle{\sum_{0< s\leq T}\widetilde{l}_s} <+\infty$,\,$P$-a.s.\\

Indeed, we just take in conditions $(\bf{A.1})\bf{(a)}$,
$(\bf{A.2})\bf{(a)}$ and $\bf{(A.3)\bf{(a)}}$, $\eta$, $C$, $\beta$
and $l$ as follows :
$$
\begin{array}{lll}
& \eta_t(\omega) = \widetilde{\eta}_t(\omega)D_t(\omega),
\\ &  C_t(\omega) = D_t(\omega),
\\ & {\beta}_t(\omega) =\widehat{\eta}_t\,\,D_t(\omega),
\\ & l_t(\omega) = \widetilde{l}_t(\omega) D_t(\omega).
\end{array}
$$ This means that the functions $f, g$ and $h$ can have, in particular, a
general growth in $(x, y)$ and stochastic quadratic growth in $z$.
\item It is not difficult to see that if $L$ or $U$ are
semimartingales, then assumption $(\bf{A.4})$ holds. Moreover if the
barriers processes $L$ and $U$ are completely separated on $[0,T[$,
\it{i.e.} $\forall t\in [0,T[$,\,\,$L_t < U_t$ and $L_{t-} < U_{t-}$
(this is equivalent to $\displaystyle{\inf_{0\leq
t<T}}(U_t-L_t)>0$), then assumption $(\bf{A.4})$ holds. Indeed, let
$$
\begin{array}{lll}
&\beta_t = \displaystyle\sup_{s\leq t} (\mid L_s\mid +\mid U_s\mid )
\\ & L_t' = \dfrac{L_t}{\beta_t}  1_{\{t<T\}} + \dfrac{L_{T-}}{\beta_{T-}}  1_{\{t=T\}}
\\ & U_t' = \dfrac{U_t}{\beta_t}  1_{\{t<T\}} + \dfrac{U_{T-}}{\beta_{T-}}
1_{\{t=T\}}.
\end{array}
$$
Then, $\forall t\in [0,T],\,\, -1\leq L_t'< U_t'\leq 1$ and $
L'_{t-}< U'_{t-}$. It follows then from the work \cite{HHO} that
there exists a semimartingale $\overline{S}$ such that $L_t'\leq
\overline{S}_t\leq U_t',\,\, \forall t\in [0,T]$. Hence, the
semimartingale $\overline{S_t}\beta_t 1_{\{t<T\}} +\xi 1_{\{t=T\}}$
is between $L_t$ and $U_t$.
\end{enumerate}
\end{remark}
\subsection{The main result}
 The following theorem constitute the main result of the paper.
 \begin{theorem}\label{thee1}
If assumptions $(\bf{A.1})$--$(\bf{A.4})$ hold then the GRBSDE
(\ref{eq0}) has a maximal solution (resp. minimal solution).
\end{theorem}
The rest of the paper is devoted to the proof of our main result
(Theorem \ref{thee1}). By means of an exponential change, the proof
of our main result consists in establishing first a correspondence
between our GRBSDE and another GRBSDE whose coefficients are more
tractable. We show that the existence of solutions for our initial
GRBSDE is equivalent to the existence of solutions for the auxiliary
GRBSDE. Since the integrability conditions on parameters are weaker,
we make use of approximations and truncations to establish the
existence result for the auxiliary GRBSDE. The final step consists
in justifying the passage to the limit and identifying the limit as
the solution of the auxiliary GRBSDE. A useful tool in our
considerations is the comparison theorem (see Theorem \ref{th111} in
Appendix). Let us start by giving equivalent forms of our GRBSDE.
%

\section{Equivalent forms of GRBSDE (\ref{eq0})}
\subsection{First equivalent form of GRBSDE (\ref{eq0})}
The following lemma shows that the existence of solutions for our
initial GRBSDE (\ref{eq0}) is equivalent to the existence of
solutions for another GRBSDE.
\begin{lemma}\label{lem1}
$(Y, Z, K^+, K^-)$ is a solution of GRBSDE (\ref{eq0}) if and only
if $(Y, Z, K^{+c}, K^{-c})$ is a solution of the following GRBSDE
\begin{equation}
\label{eq1} \left\{
\begin{array}{ll}
(i) & 
 Y_{t}=\xi
+\integ{t}{T}f(s,Y_{s},Z_{s})ds+\dint_t^Tg(s,
Y_s)dA_s -\sum_{t<s\leq T}\Delta Y_s\\
&\qquad\quad+\integ{t}{T}dK_{s}^{+c} -\integ{t}{T}dK_{s}^{-c}
-\integ{t}{T}Z_{s}dB_{s}\,, t\leq T,
\\ (ii)&
\forall t\in]0,T],\,\, Y_{t-} = L_{t-} \vee [Y_t + h(t, Y_{t-}, Y_t)]\wedge U_{t-},\\
(iii)& \integ{0}{T}( Y_{t}-L_{t})
dK_{t}^{+c}= \integ{0}{T}( U_{t}-Y_{t}) dK_{t}^{-c}=0,\,\, \mbox{a.s.}, \\
(iv)& Y\in {\cal D}, \quad K^{+c}, K^{-c}\in {\cal K}^c, \quad Z\in
{\cal L}^{2,d},  \\ (v)& dK^{+c}\perp  dK^{-c}.
\end{array}
\right. \end{equation}
\end{lemma}
\bop. Suppose that $(Y, Z, K^+, K^-)$ is a solution of GRBSDE
(\ref{eq0}) and let $K^{\pm c}$ be the continuous part of $K^{\pm}$.
Clearly for every $t\in [0, T]$
$$
Y_{t}=\xi +\integ{t}{T}f(s,Y_{s},Z_{s})ds+\dint_t^Tg(s, Y_s)dA_s
-\sum_{t<s\leq T}\Delta_s Y+\integ{t}{T}dK_{s}^{+c}
-\integ{t}{T}dK_{s}^{-c} -\integ{t}{T}Z_{s}dB_{s},
$$
with $\Delta_t Y = -h(t, Y_{t-}, Y_t )+\Delta_t K^- -\Delta_t K^+$.
Then for each $t\in ]0, T]$,\, $$ Y_{t-} =Y_t+h(t, Y_{t-}, Y_t
)+\Delta_t K^+ -\Delta_t K^-.$$ We should remark first that
\begin{equation}\label{eqq13}
 dK^+ \perp dK^- \Longleftrightarrow dK^{+c} \perp
 dK^{-c}
 \,\,\,\mbox{and}\,\,\, \Delta K^- \Delta K^+ =0.
 \end{equation}
 Now, we distinguish the following three cases
:
\begin{enumerate}
\item If $\Delta_t K^+ >0$,
since $dK^+ \perp dK^-$ and $\integ{0}{T}( Y_{t-}-L_{t-}) dK_{t}^+=
0$ it follows that $\Delta_t K^- = 0$ and $Y_{t-} = L_{t-}$.
Consequently $\Delta_t K^+  = \bigg(L_{t-} - [Y_t+h(t, Y_{t-}, Y_t
)]\bigg)^+ $ and $L_{t-}
> Y_t+h(t, Y_{t-}, Y_t )$. Henceforth for each $t\in ]0, T]$,
$$
Y_{t-} = L_{t-} \vee [Y_t + h(t, Y_{t-}, Y_t)]\wedge U_{t-}.
$$
\item
If $\Delta_t K^- >0$, by the same way as above we get $\Delta_t K^-
= \bigg(U_{t-}- [Y_t+h(t, Y_{t-}, Y_t )]\bigg)^- $ and then for each
$t\in ]0, T]$
$$
Y_{t-} = L_{t-} \vee [Y_t + h(t, Y_{t-}, Y_t)]\wedge U_{t-}.
$$
\item If $\Delta_t K^- = \Delta_t K^+=0$, we have $Y_{t-} = Y_t +
h(t, Y_{t-}, Y_t)\in [L_{t-}, U_{t-}]$ and then for each $t\in ]0,
T]$
$$
Y_{t-} = L_{t-} \vee [Y_t + h(t, Y_{t-}, Y_t)]\wedge U_{t-}.
$$
\end{enumerate}
Consequently in the three cases we obtain for each $t\in ]0, T]$
\begin{equation}\label{equa1}
\begin{array}{ll}
 & \Delta_t K^+  =
\bigg(L_{t-}- [Y_t+h(t, Y_{t-}, Y_t )]\bigg)^+
\\ & \Delta_t K^-  =
\bigg(U_{t-}- [Y_t+h(t, Y_{t-}, Y_t )]\bigg)^-
\\
& Y_{t-} = L_{t-} \vee [Y_t + h(t, Y_{t-}, Y_t)]\wedge U_{t-}.
 \end{array}
\end{equation}
 Hence $(ii)$ of Equation (\ref{eq1}) is satisfied. \\ Moreover, since
 $
 0\leq (Y_t - L_t) dK^{+c}_t\leq  (Y_{t-} - L_{t-}) dK^{+}_t =
 0,\,\,\mbox{and}\,\, 0\leq (U_t - Y_t) dK^{-c}_t\leq  (U_{t-} - Y_{t-})
 dK^{-}_t =0,
 $
 then we get $(iii)$ of Equation (\ref{eq1}). \\ In view of (\ref{eqq13}), we have also $(v)$ of Equation
 (\ref{eq1}). Hence $(Y, Z, K^{+c}, K^{-c})$ is a solution of GRBSDE
 (\ref{eq1}).

 On another hand, suppose now that $(Y, Z, K^{+c}, K^{-c})$ is a solution
 of GRBSDE (\ref{eq1}) and set for each
$t\in]0, T]$
\begin{equation}\label{eq11}
\begin{array}{ll}
 & \Delta_t K^+  =
\bigg(L_{t-}- [Y_t+h(t, Y_{t-}, Y_t )]\bigg)^+
\\ & \Delta_t K^-  =\bigg(U_{t-}- [Y_t+h(t, Y_{t-},
Y_t )]\bigg)^-.
 \end{array}
 \end{equation} By $(ii)$ of Equation (\ref{eq1}) we have for each $t\in]0, T]$
$$
\begin{array}{ll}
  -\Delta_t Y & =Y_{t-} -Y_t = (L_{t-} -Y_t)\vee h(t, Y_{t-},
Y_t )\wedge (U_{t-} -Y_t),
\\ & = h(t, Y_{t-}, Y_t) +[(U_{t-} -Y_t- h(t, Y_{t-}, Y_t))\wedge
0\vee (L_{t-} -Y_t- h(t, Y_{t-}, Y_t))]
\\ & = h(t, Y_{t-}, Y_t) + \Delta_t K^+ - \Delta_t K^-.
 \end{array}
 $$
 Hence
$$
Y_{t}=\xi +\integ{t}{T}f(s,Y_{s},Z_{s})ds+\dint_t^Tg(s, Y_s)dA_s
+\sum_{t<s\leq T}h(s, Y_{s-}, Y_s)+\integ{t}{T}dK_{s}^+
-\integ{t}{T}dK_{s}^- -\integ{t}{T}Z_{s}dB_{s},
$$
where $K^{\pm}$ is defined by : $K^{\pm}_t= K_t^{\pm
c}+\displaystyle{\sum_{0<s\leq t}\Delta_sK^{\pm}}$, $\forall t\in
[0, T]$. Moreover it follows that
\begin{enumerate}
\item For each $t\in[0, T]$, $L_t\leq Y_t\leq U_t$ since $L_{t-}\leq
Y_{t-}\leq U_{t-}$ and $Y_T=L_T=U_T=\xi$.
\item $\Delta_t K^+\Delta_t K^- = 0$ since
$L_{t-}\leq U_{t-}$ and then $dK^+ \perp dK^-$.
\item $\integ{0}{T}( Y_{t-}-L_{t-}) dK_{t}^{+}=\integ{0}{T}(
Y_{t}-L_{t}) dK_{t}^{+c} + \displaystyle{\sum_{0< s\leq t}}(
Y_{s-}-L_{s-})\Delta_s K^+  = \displaystyle{\sum_{0< s\leq t}}(
Y_{s-}-L_{s-})\Delta_s K^+ =0$,  since if $Y_{s-} > L_{s-}$ then
$Y_s + h(s, Y_{s-}, Y_s) > L_{s^{-}}$ and hence $\Delta_s K^+=0$.
Similarly it follows also that $\integ{0}{T}( U_{t-}-Y_{t-})
dK_{t}^{-}= 0$. Therefore $(Y, Z, K^+, K^-)$ is a solution of GRBSDE
(\ref{eq0}).
\end{enumerate}
 \eop
\begin{remark}\label{lem2}
The maximal solution $(Y, Z, K^{+c}, K^{-c})$ of GRBSDE (\ref{eq1})
satisfies, $P-a.s.$\, $\forall t\in ]0, T]$
$$
Y_{t- } = \max\{ x\in[L_{t-}, U_{t-}]\,\, : \,\, x= L_{t-}\vee[ Y_t+
h(t, x, Y_t)]\wedge U_{t-} \}.
$$
Indeed, let $\tau\in ]0,T]$ be a stopping time and set
$$
\overline{\xi} = \max\{ x\in[L_{\tau -}, U_{\tau -}]\,\, : \,\, x=
L_{\tau -}\vee [Y_\tau + h(\tau, x, Y_\tau)]\wedge U_{\tau -} \}.
$$
Let $(Y', Z', K'^{+c}, K'^{-c})$ be a solution of GRBSDE
(\ref{eq1})(which is exists according to our main result) associated
to the data : $f'= 1_{\{s\leq \tau\}} f$,\,\, $g'= 1_{\{s\leq
\tau\}} g$,\,\, $h'= 1_{\{s< \tau\}} h$,\,\, $\xi =
\overline{\xi}$,\,\, $L'_t = L_t 1_{\{t< \tau\}} + L_{\tau
-}1_{\{t\geq \tau\}}$ and $U'_t = U_t 1_{\{t< \tau\}} + U_{\tau
-}1_{\{t\geq \tau\}}$. Set also
$$
\begin{array}{ll}
 & Y''_t = Y'_t 1_{\{t<
\tau\}} + Y_{t}1_{\{t\geq \tau\}}
\\ &  Z''_t = Z'_t 1_{\{t<
\tau\}} + Z_{t}1_{\{t\geq \tau\}}
\\ &  dK''^{\pm c}_t = dK'^{\pm c}_t 1_{\{t<
\tau\}} + dK_{t}^{\pm c}1_{\{t\geq \tau\}}.
 \end{array}
 $$
Clearly $(Y'', Z'', K''^{+c}, K''^{-c})$ is also a solution of
GRBSDE (\ref{eq1}). Since $(Y, Z, K^{+c}, K^{-c})$ is a maximal
solution of GRBSDE (\ref{eq1}), then $Y''_t \leq Y_t, \,\,\forall
t\in [0, T]$. Henceforth
$$
\begin{array}{ll}
Y''_{\tau -} &=Y'_{\tau -} = \overline{\xi} \\ & \leq Y_{\tau -} =
L_{\tau -}\vee (Y_{\tau}+h(\tau, Y_{\tau -}, Y_\tau))\wedge U_{\tau -}\\
& \leq \max\{ x\in[L_{\tau -}, U_{\tau -}]\,\, : \,\, x= L_{\tau
-}\vee(Y_{\tau}+ h(\tau, x, Y_\tau))\wedge U_{\tau -}
\}=\overline{\xi}.
\end{array}
$$
Then, for every stopping time $\tau\in ]0, T]$ we get
$$
Y_{\tau -}=\max\{ x\in[L_{\tau -}, U_{\tau -}]\,\, : \,\, x= L_{\tau
-}\vee (Y_{\tau}+h(\tau, x, Y_\tau))\wedge U_{\tau -} \},\,\,\,
P-a.s.
$$
It therefore follows that, $P-a.s.$\, $\forall t\in ]0, T]$
$$
Y_{t- } = \max\{ x\in[L_{t-}, U_{t-}]\,\, : \,\, x= L_{t-}\vee[ Y_t+
h(t, x, Y_t)]\wedge U_{t-} \}.
$$
\end{remark}

\subsection{Second equivalent form of GRBSDE (\ref{eq0})}\label{sec3.2} In this part,
by using an exponential transform, we transform the GRBSDE with two
obstacles into another equivalent one whose data satisfy some "good"
conditions. This transformation allow us, in
particular, to bound the terminal condition and the barriers associated with the transformed GRBSDE. \\
To begin with, let $m\in {\cal K} +\R_+$ and suppose that GRBSDE
(\ref{eq0}) has a solution. It follows then from It\^{o}'s formula
that
$$
\begin{array}{lll}
 (Y_t-S_t-m_t)m_t & = (\xi-S_T-m_T)m_T + \dint_{t}^{T}m_s
f(s,Y_{s},Z_{s})ds+ \dint_{t}^{T}m_s g(s,Y_{s})dA_s \\ &+
\dint_{t}^{T} m_sdK_s^{+c} -\dint_{t}^{T}m_s dK_s^{-c}
-\dint_{t}^{T}m_s(Z_s-\gamma_s) dB_s + \dint_{t}^{T}m_s dV^c_s\\
&+\dint_{t}^{T}m_s dm_s^c-\dint_{t}^{T}(Y_s-S_s-m_s)dm_s^c
-\displaystyle{\sum_{t<s\leq T}}\Delta_s [(Y_.-S_.-m_.)m_.].
\end{array}
 $$
Setting $e_t:= e^{m_t(Y_t-S_t-m_t)}$, it follows that
$$
\begin{array}{lll}
& e_t:= e^{(Y_t-S_t-m_t)m_t} \\ &= e_T+ \dint_{t}^{T}e_sm_s
f(s,Y_{s},Z_{s})ds+ \dint_{t}^{T}e_s m_s g(s,Y_{s})dA_s +
\dint_{t}^{T} e_sm_sdK_s^{+c} -\dint_{t}^{T}e_s m_s dK_s^{-c}\\ &
-\dint_{t}^{T}e_sm_s(Z_s-\gamma_s) dB_s + \dint_{t}^{T}e_sm_s
dV^c_s+\dint_{t}^{T}e_sm_s dm_s^c-\dint_{t}^{T}e_s
(Y_s-S_s-m_s)dm_s^c \\
&-\dfrac12\dint_{t}^{T}e_sm_s^2|Z_s-\gamma_s|^2 ds
-\displaystyle{\sum_{t<s\leq T}}\Delta_s e_..
\end{array}
 $$
 Since $d(e^{m_t^2})^c = 2m_s e^{-m_s^2} dm_s^c$, then
\begin{equation}\label{eq6}
\begin{array}{lll}
&e_t -e^{-m_t^2} \\ & = e_T-e^{-m_T^2}+ \dint_{t}^{T}e_sm_s\bigg(
f(s,Y_{s},Z_{s})-\frac12 m_s|Z_s-\gamma_s|^2 \bigg) ds
\\ &+
\dint_{t}^{T}\bigg(e_s m_s g(s,Y_{s})dA_s+ e_sm_s dV^c_s +e_sm_s
dm_s^c -e_s (Y_s-S_s-m_s)dm_s^c-2m_s e^{-m_s^2} dm_s^c\bigg)
\\& + \dint_{t}^{T} e_sm_sdK_s^{+c}-\dint_{t}^{T}e_s m_s dK_s^{-c} -\dint_{t}^{T}e_sm_s(Z_s-\gamma_s)
dB_s -\displaystyle{\sum_{t<s\leq T}}\Delta_s (e_. -e^{-m_.^2}).
\end{array}
\end{equation}

\ni Let $|V|$ denotes the total variation of the process $V$ and
choose the process $m$ as follows : $\forall s\in [0, T]$,
$$
\begin{array}{lll}
&\ni m_s = 4\bigg[\displaystyle\sup_{r\leq
s}\bigg(|U_r|+|C_r|+|L_r|\bigg)+|V|_s+\dint_0^s(1+\eta_r+|\gamma_r|^2)dr
+\dint_0^s(1+\beta_r)dA_r +\displaystyle\sum_{0<r\leq s}
l_r+1\bigg].
\end{array}
$$
Define also for every $s\in [0, T]$,
\begin{equation}\label{equa2}
\begin{array}{lll}
  &\ni \bullet\quad\overline{\xi}=e_T -e^{-m_T^2} =e^{
 m_T(\xi-S_T-m_T)}-e^{-m_T^2}. 
 \\ &\ni \bullet\quad
\overline{L}_s = e^{ m_s(L_s-S_s-m_s)}-e^{-m_s^2}.\\ &\ni
\bullet\quad \overline{U}_s = e^{ m_s(U_s-S_s-m_s)}-e^{-m_s^2}. \\
&\ni \bullet\quad \overline{Y}_s = e_s -e^{-m_s^2} =
e^{m_s(Y_s-S_s-m_s)}-e^{-m_s^2}.\\ &\ni \bullet\quad \overline{Z}_s
= m_s (\overline{Y}_s +e^{-m_s^2})(Z_s-\gamma_s).
\\ &\ni \bullet\quad d\overline{K}_s^{\pm c} = m_s (\overline{Y}_s
+e^{-m_s^2})dK^{\pm c}_s.
\\ &  \ni \bullet\quad \Delta_s \overline{K}^+ = \bigg(\overline{L}_{s-}-
e^{[(Y_{s}+h(s,Y_{s-} ,Y_s)-S_{s-}-m_{s-})m_{s-}]} +
e^{-m_{s-}^2}\bigg)^+.
\\ &  \ni \bullet\quad \Delta_s \overline{K}^{\,-} = \bigg(\overline{U}_{s-}-
e^{[(Y_{s}+h(s,Y_{s-} ,Y_s)-S_{s-}-m_{s-})m_{s-}]} +
e^{-m_{s-}^2}\bigg)^-.
\\ &  \ni \bullet\quad \overline{K}_s^{\pm} = \overline{K}_s^{\pm c}
+\displaystyle{\sum_{0 <r\leq s}}\Delta_r\overline{ K}^{\pm}.
\end{array}
\end{equation}
\begin{remark}
\begin{enumerate}
\item It should be noted that $m$ is ${\cal
F}_t$-adapted, \it{rcll} and increasing process. 
\item Since $S_T=\xi$ then $\overline{\xi}= 0$.
\item Since, $\forall s\in ]0, T]$,\,\,  $Y_{s-} = L_{s-}\vee [Y_s+
h(s,Y_{s-} ,Y_s)]\wedge U_{s-}$,
then for every $s\in]0, T]$,
$$
\begin{array}{lll}
 \overline{Y}_{s-} &= \overline{L}_{s-}\vee\bigg[
e^{[(Y_{s}+h(s,Y_{s-} ,Y_s)-S_{s-}-m_{s-})m_{s-}]} -
e^{-m_{s-}^2}\bigg]\wedge \overline{U}_{s-}.
\end{array}
$$
\end{enumerate}
\end{remark}
\ni Coming back to Equation (\ref{eq6}), it is clear that the GRBSDE
(\ref{eq0}) can be written as follows :
\begin{equation}
\label{eq51} \left\{
\begin{array}{ll}
(i) & 
 \overline{Y}_{t}=\integ{t}{T}\widetilde{f}(s,\overline{Y}_{s},\overline{Z}_{s})ds +
 \integ{t}{T}\widetilde{g}(s,\overline{Y}_{s})d\overline{A}_s -\displaystyle{\sum_{t<s\leq
T}}\Delta_s \overline{Y}  +\integ{t}{T}d\overline{K}_{s}^{+c}
\\
&\qquad\quad-\integ{t}{T}d\overline{K}_{s}^{\,-c}
-\integ{t}{T}\overline{Z}_{s}dB_{s} \,, t\leq T,
\\ (ii)&\forall s\in]0,T],\,\, \overline{Y}_{s-} = \overline{L}_{s-}\vee\bigg[
\overline{Y}_{s}+ \widetilde{h}(s, \overline{Y}_{s-},
\overline{Y}_{s})\bigg]\wedge \overline{U}_{s-}
 \\ (iii)& \integ{0}{T}( \overline{Y}_{t}-\overline{L}_{t})
d\overline{K}_{t}^{+c}= \integ{0}{T}(
\overline{U}_{t}-\overline{Y}_{t}) d\overline{K}_{t}^{-c}=0,\,\,
\mbox{a.s.} \\ (iv)& \overline{Y}\in {\cal D}, \quad
\overline{K}^{\,+c}, \overline{K}^{\,-c}\in {\cal K}^c, \quad
\overline{Z}\in {\cal L}^{2,d},
 \\ (v)&
d\overline{K}^+\perp d\overline{K}^-,
\end{array}
\right. \end{equation} where $\overline{A}$,
 $\widetilde{g}$, $\widetilde{f}$ and $\widetilde{h}$ are given by
 : for each $ s\in [0, T], y\geq\overline{L}_s$,
$x\geq\overline{L}_{s-}$ and $z\in\R^d$
$$
\begin{array}{lll}
  &\bullet \quad \overline{A}_s =  2 \dint_0^s e^{-m_r} dm_r^c,
\\
& \\ &\bullet\quad \widetilde{f}(s,y,z)= \\ & \quad m_s(y+
e^{-m_s^2}) \bigg(f(s, \dfrac{\ln(y+e^{-m_s^2})}{m_s}+m_s +S_s,
\dfrac{z}{m_s (y+e^{-m_s^2})}+\gamma_s)\bigg)
-\dfrac{|z|^2}{2(y+e^{-m_s^2})},
\\ & \\ & \bullet \quad \widetilde{g}(s,y) =
m_s(y+e^{-m_s^2}) \bigg[
g(s,\dfrac{\ln(y+e^{-m_s^2})}{m_s}+m_s+S_s)\dfrac{dA_s}{d\overline{A}_s}
+(\dfrac{dV^c_s}{d\overline{A}_s}+\dfrac{dm^c_s}{d\overline{A}_s})\bigg]\\
& \qquad\qquad\qquad
-2m_se^{-m_s^2}\dfrac{dm^c_s}{d\overline{A}_s}-(\overline{y}+e^{-m_s^2})
\dfrac{\ln(y+e^{-m_s^2})}{m_s}\dfrac{dm^c_s}{d\overline{A}_s},
\\ & \\ &\bullet
\quad \widetilde{h}(s,x, y)  = \overline{\overline{h}}(s,x,y) -y
-e^{-m_{s-}^2},\quad \mbox{where}\\ &
\overline{\overline{h}}(s,x,y)= e^{\bigg[\dfrac{\ln(y
+e^{-m_s^2})}{m_s} + h\bigg(s, \dfrac{\ln(x +e^{-m_{s-}^2})}{m_{s-}}
+S_{s-} +m_{s-},\dfrac{\ln(y +e^{-m_s^2})}{m_s} +S_s +m_s\bigg)
+\Delta_s S +\Delta_s m\bigg]m_{s-}}.
\end{array}
$$
 It follows from Equation (\ref{eq51}) and Lemma \ref{lem1} that the GRBSDE (\ref{eq0})
can be written also as follows :
\begin{equation}
\label{eq5} \left\{
\begin{array}{ll}
(i) & 
 \overline{Y}_{t}=\integ{t}{T}\overline{f}(s,\overline{Y}_{s},\overline{Z}_{s})ds +
 \integ{t}{T}\overline{g}(s,\overline{Y}_{s})d\overline{A}_s +\displaystyle{\sum_{t<s\leq
T}}\overline{h}(s,\overline{Y}_{s-} ,\overline{Y}_s) + \integ{t}{T} d\overline{R}_s\\
&\qquad\quad+\integ{t}{T}d\overline{K}_{s}^+
-\integ{t}{T}d\overline{K}_{s}^{\,-}
-\integ{t}{T}\overline{Z}_{s}dB_{s} \,, t\leq T,
\\ (ii)& \forall t\leq T,\,\, \overline{L}_{t}\leq \overline{Y}_{t}\leq \overline{U}_{t}
 \\ (iii)& \integ{0}{T}( \overline{Y}_{t-}-\overline{L}_{t-})
d\overline{K}_{t}^+= \integ{0}{T}(
\overline{U}_{t-}-\overline{Y}_{t-}) d\overline{K}_{t}^-=0,\,\,
\mbox{a.s.} \\ (iv)& \overline{Y}\in {\cal D}, \quad
\overline{K}^{\,+}, \overline{K}^{\,-}\in {\cal K}, \quad
\overline{Z}\in {\cal L}^{2,d},
 \\ (v)&
d\overline{K}^+\perp d\overline{K}^-,
\end{array}
\right. \end{equation} where $\overline{f}$, $\overline{g}$,
$\overline{h}$ and $\overline{R}$ is given by
 : for each $ s\in [0, T]$, $\overline{y}\in\R$, $\overline{x}\in\R$
 and $\overline{z}\in\R^d$
$$
\begin{array}{lll}
   &\bullet\quad \overline{f}(s,\overline{y},\overline{z})
=\widetilde{f}(s,\overline{L}_s\vee
\overline{y}\wedge\overline{U}_s,\overline{z}) -\dfrac12
\overline{\eta}_s, \quad\mbox{where}\quad \overline{\eta}_s=
2e^{-m_s}( \eta_s+ \gamma_s^2),
\\ & \\ & \bullet \quad \overline{g}(s,\overline{y}) ={\widetilde{g}(s,
\overline{L}_s\vee \overline{y}\wedge\overline{U}_s)-\dfrac12},
\\ & \\ &\bullet\quad \overline{h}(s, \overline{x}, \overline{y}) =
\widetilde{h}(s, \overline{L}_{s-}\vee
\overline{x}\wedge\overline{U}_{s-}, \overline{L}_s\vee
\overline{y}\wedge\overline{U}_s)-
 e^{3m_s^2}\Delta_s m,
\\ & \\ & \bullet\quad  \overline{R}_s := \dfrac12 \dint_0^sd\overline{A}_r +
\dfrac12 \dint_0^s\overline{\eta}_r dr + \displaystyle{\sum_{0<r\leq
s}} e^{3m_r^2}\Delta_r m.
\end{array}
$$
\begin{remark}\label{rem0}
It is clear now that if $(Y, Z, K^+, K^-)$  is a solution (resp.
maximal solution) of GRBSDE (\ref{eq0}) associated with the data
$(\xi, f, g, A, L, U)$ then $(\overline{Y}, \overline{Z},
\overline{K}^+, \overline{K}^-)$, defined by (\ref{equa2}) and
associated with coefficient the data $(\overline{\xi}, \overline{f},
\overline{g}, \overline{h}, \overline{A}, \overline{R},
\overline{L}, \overline{U})$, is a solution (resp. maximal solution)
of GRBSDE $(\ref{eq5})$. Conversely, Suppose that there exists a
solution (resp. maximal solution)
$(\overline{Y},\overline{Z},\overline{K}^+,\overline{K}^-)$ for
GRBSDE $(\ref{eq5})$. Hence, by setting, for all $t\leq T$
$$
\begin{array}{ll}
&Y_t = \dfrac{\ln(\overline{Y}_t+e^{-m_t^2})}{m_t} +S_t+m_t,
\\ &
Z_t = \dfrac{\overline{Z}_t}{m_t ({\overline{Y}_t+e^{-m_t^2}  })}
+\gamma_t,
\\ &
 K_t^{\pm} = \dint_0^t\dfrac{d\overline{K}_s^{\pm c}}{m_s(
{\overline{Y}_s}+e^{-m_s^2})} + \displaystyle{\sum_{0 <r\leq
t}}\Delta_r K^{\pm},
\end{array}
$$
where $\Delta_t K^+  = \bigg(L_{t-}- [Y_t+h(t, Y_{t-}, Y_t
)]\bigg)^+$\,\, and $\Delta_t K^-  = \bigg(U_{t-}- [Y_t+h(t, Y_{t-},
Y_t )]\bigg)^-$, it is clear that $(Y,Z,K^+,K^-)$ is a solution
(resp. maximal solution) for GRBSDE $(\ref{eq0})$.
\end{remark}
The following proposition states some properties on the data
$(\overline{\xi}, \overline{f}, \overline{g}, \overline{h},
\overline{A}, \overline{R}, \overline{L}, \overline{U})$ of the
transformed GRBSDE (\ref{eq5}).
\begin{proposition}\label{pro0} Assume that assumptions
$(\bf{A.1})-(\bf{A.3})$ hold. Then we have the following :
\begin{enumerate}
\item The function $\overline{f}$ is $\cal P$-measurable and
continuous with respect to $(\overline{y},\overline{z})$ satisfying
for every $s\in [0,T],\,\,$ $\overline{y}\in \R$ and
$\overline{z}\in \R^d $
\begin{equation}\label{E2} -\overline{{\eta}}_s -e^{2m_s^2}|\overline{z}|^2\leq \overline{f}(s,
\overline{y},\overline{z}) \leq 0.
\end{equation}
\item 
$\dint_0^T \overline{\eta}_s ds\leq  \overline{A}_T\leq 1$.
\item The function $\overline{g}$ is $\cal P$-measurable and
continuous with respect to $y$ satisfying, for all $s\in [0,T]$ and
$\overline{y}\in \R,$
\begin{equation}\label{eq19}
-1\leq \overline{g}(s, \overline{y})\leq 0.
\end{equation}
\item For all $t\in [0,T],\,\,\, -1\leq \overline{L}_t \leq 0 \leq\overline{U}_t \leq 1$.
\item The function $\overline{h}$ satisfies the following
properties
\begin{enumerate}
\item $P$-a.s.,
$\forall (t, x)\in [0, T]\times\R$, the function
$\overline{y}\mapsto \overline{y}+\overline{h}(t,
\overline{x},\overline{y}) $ is nondecreasing and continuous on
$\R$.
\item $P$-a.s.,
$\forall t\in[0, T]$,\, $\forall \overline{y}\in\R$,  the function
$\overline{x}\mapsto \overline{h}(t, \overline{x},\overline{y}) $ is
continuous on $\R$.
\item $P$-a.s.,
$\forall t\in[0, T]$,\, $\forall
\overline{x}\in\R,\,\forall\overline{y}\in\R$,
$$
-e^{3m_t^2}\Delta_t m +\Delta_t e^{-m^2_.}\leq
\overline{h}(t,\overline{x}, {\overline{y}})\leq 0.
$$
\end{enumerate}
 \end{enumerate}
 \end{proposition}
\bop.  $1.$ It is not difficult to see that $\overline{f}$ is $\cal
P$-measurable and continuous with respect to $(y,z)$ on $\R\times
\R^d$ for every $(t,\omega)\in ]0, T[\times\Omega$, since $f$ is.
Let us prove Inequality (\ref{E2}). Let $s\in [0,T],
\,\,\overline{y}\in [\overline{L}_s(\omega),
\overline{U}_s(\omega)]$ and $\overline{z}\in \R^d$. By condition
$(\bf{A.1})$ we get
$$\begin{array}{ll}\widetilde{f}(s,\omega, \overline{y},\overline{z})
& \leq m_s(\overline{y}+e^{-m_s^2})\bigg( \eta_s+ C_s
\bigg|\dfrac{|\overline{z}|}{m_s(\overline{y}+e^{-m_s^2})}+\gamma_s\bigg|^2\bigg)
-\dfrac{|\overline{z}|^2}{2(\overline{y}+e^{-m_s^2})}\\
&\leq m_s(\overline{y}+e^{-m_s^2})( \eta_s+ 2C_s\gamma_s^2) +2C_s
\dfrac{|\overline{z}|^2}{m_s(\overline{y}+e^{-m_s^2})}-\dfrac{|\overline{z}|^2}{2(\overline{y}+e^{-m_s^2})} \\
& \leq m_s(\overline{U}_s+e^{-m_s^2})( \eta_s+ 2C_s\gamma_s^2)
-(m_s-4C_s) \dfrac{|\overline{z}|^2}{2m_s(\overline{y}+e^{-m_s^2})}
\\ & \leq m_s e^{ m_s(U_s-S_s-m_s)}( \eta_s+ 2C_s\gamma_s^2)
\\ & \leq \frac{m_s^2}{2}e^{-\frac{m_s^2}{2}}( \eta_s+ \gamma_s^2)
\\ &\leq 2\,\displaystyle{\sup_{x\geq 4}(xe^{-x})}e^{-m_s}( \eta_s+ \gamma_s^2)
\\ & \leq e^{-m_s}( \eta_s+ \gamma_s^2) = \dfrac{\overline{\eta}_s}{2}.
\end{array}$$
where we have used the elementary inequality $(a+ b)^2 \leq 2(a^2+
b^2)$ and the fact that $m_s-4C_s \geq 0$ and $U_s -S_s-m_s\leq
-\frac{m_s}{2}$, $\quad \forall s\in [0,T]$ .\\ On the other hand,
by using condition $(\bf{A.1})$, we get also that
$$\begin{array}{ll} \widetilde{f}(s,\omega, \overline{y},\overline{z})
& \geq m_s(\overline{y}+e^{-m_s^2})\bigg( -\eta_s- C_s
\bigg|\dfrac{|\overline{z}|}{m_s(\overline{y}+e^{-m_s^2})}
+\gamma_s\bigg|^2\bigg)-\dfrac{|\overline{z}|^2}{2(\overline{y}+e^{-m_s^2})}
\\ & \geq -\frac12\overline{\eta}_s - 2
m_s(\overline{y}+e^{-m_s^2})C_s\dfrac{|\overline{z}|^2}{m_s^2(\overline{y}+e^{-m_s^2})^2}
-\dfrac{|\overline{z}|^2}{2(\overline{y}+e^{-m_s^2})}
\\ & \geq -\frac12\overline{\eta}_s -(\dfrac{2C_s}{m_s}+\dfrac12) \dfrac{|\overline{z}|^2}{(\overline{y}+e^{-m_s^2})}
\\ & \geq -\frac12\overline{\eta}_s -\dfrac{|\overline{z}|^2}{(\overline{y}+e^{-m_s^2})}
 \\ & = -\frac12\overline{\eta}_s - e^{ -m_s(L_s-S_s-m_s)} |\overline{z}|^2
 \\ & \geq -\frac12\overline{\eta}_s - e^{ \frac{3}{2}m_s^2} |\overline{z}|^2
\end{array}
$$
since $\overline{y} +e^{ m_s^2}\geq \overline{L}_s+e^{ m_s^2} = e^{
-m_s(L_s-S_s-m_s)}$ and $\frac{2C_s}{m_s}+\frac{1}{2} \leq 1$,$\quad
\forall s\in [0,T]$. Consequently
$$ -\overline{{\eta}}_s
-e^{2m_s^2}|\overline{z}|^2\leq \overline{f}(s,
\overline{y},\overline{z}) \leq 0.
$$
2. Since for every $s\in [0,T]$,\,\, $\dfrac{(\eta_s
+\gamma_s^2)ds}{dm_s^c}\leq \frac14$, then we have
$$
\dint_0^T \overline{\eta}_s ds =2 \dint_0^T e^{-m_s}\dfrac{(\eta_s
+\gamma_s^2) ds}{dm_s^c} dm_s^c
\leq  \dfrac{1}{2} \dint_0^T e^{-m_s}dm_s^c.
$$
Since $ \dint_0^T e^{-m_s} dm_s^c =e^{-m_0}
-e^{-m_T}+\displaystyle{\underbrace{\sum_{0\leq s\leq T} \Delta
e^{-m_s}}_{\leq 0} \leq  e^{-m_0}}$, then $ \dint_0^T
\overline{\eta}_s ds\leq \dfrac{1}{4}\overline{A}_T\leq
e^{-m_0}\leq e^{-1}\leq 1. $\\
3. By using a similar calculation as above it is easy to prove that
for all $s\in [0,T]$ and $\overline{y}\in [\overline{L}_s(\omega),
\overline{U}_s(\omega)]$
$$
-1\leq \overline{g}(s, \overline{y})\leq 0.
$$
4. It is not difficult to prove that  $$\forall t\in [0,T],\,\,\,
-1\leq \overline{L}_t \leq 0 \leq\overline{U}_t \leq 1.$$
 5. $(a)$ and $(b)$ are obvious. Let us prove $(c)$. By
definition of process $m$ we have for all $s\in ]0,T]$,
$\overline{x}\in [\overline{L}_{s-}(\omega),
\overline{U}_{s-}(\omega)]$ and $\overline{y}\in
[\overline{L}_s(\omega), \overline{U}_s(\omega)]$
$$\begin{array}{ll}
&\widetilde{h}(s,\overline{x}, \overline{y})
\\ &  \geq
e^{\bigg[\dfrac{\ln(\overline{y} +e^{-m_s^2})}{m_s} - l_s +\Delta_s
S +\Delta_s m\bigg]m_{s-}} -(\overline{y} +e^{-m_s^2})+ \Delta_s
e^{-m_.^2}
\\ & = (\overline{y} +e^{-m_s^2})\bigg(e^{\bigg[-l_s +\Delta_s V +\Delta_s m\bigg]m_{s-}}e^{\frac{-\Delta_s
m\ln(\overline{y} +e^{-m_s^2})}{m_s}} -1\bigg) +\Delta_s e^{-m_.^2}
\\ & \geq \Delta_s e^{-m_.^2}.
\end{array}
$$
On the other hand, by definition of process $m$ we have for all
$s\in ]0,T]$, $\overline{x}\in [\overline{L}_{s-}(\omega),
\overline{U}_{s-}(\omega)]$ and $\overline{y}\in
[\overline{L}_s(\omega), \overline{U}_s(\omega)]$
$$
\begin{array}{lll}
 \widetilde{h}(s, \overline{x}, \overline{y})
 &\leq
e^{\bigg[\dfrac{\ln(\overline{y} +e^{-m_s^2})}{m_s} + l_s +\Delta_s
S +\Delta_s m\bigg]m_{s-}} -(\overline{y} +e^{-m_s^2})
\\ & \leq (\overline{y} +e^{-m_s^2})\bigg(e^{\bigg[l_s+\Delta_s V +\Delta_s m\bigg]m_{s-}}e^{\frac{-\Delta_s
m\ln(\overline{y} +e^{-m_s^2})}{m_s}} -1\bigg)
\\ & \leq e^{ m_s(U_s-S_s-m_s)}(e^{\frac{3}{2}(\Delta_sm)m_{s-} }\,\,e^{\frac{3}{2}(\Delta_sm)m_{s-} }-1)
\\ & \leq e^{ -\frac{m_s^2}{2}}(e^{3m_{s}\Delta_sm }-1)
\\ & \leq e^{ -\frac{m_s^2}{2}}\bigg(3m_{s}\Delta_sm e^{3m_{s}\Delta_sm } \bigg)

 \\ & \leq 3m_{s}e^{ -\frac{m_s^2}{2}}(\Delta_sm) e^{3m_{s}^2 }

\\ & \leq \displaystyle{\sup_{x\geq 4}(3xe^{-\frac{x^2}{2}}})(\Delta_sm)e^{3m_{s}^2}
\\ & \leq (\Delta_sm) e^{3m_{s}^2}
\end{array}
$$

\ni The proof of Proposition \ref{pro0} is then finished.
 \eop

\subsection{An equivalent result to Theorem \ref{thee1}}
Now, by taking advantage of the previous analysis, especially Remark
\ref{rem0} and Proposition \ref{pro0}, our problem is then reduced
to find the maximal solution of the following GRBSDE :
\begin{equation}
\label{eq221} \left\{
\begin{array}{ll}
(i) & 
 {Y}_{t}=\integ{t}{T}{f}(s,{Y}_{s},{Z}_{s})ds +
 \integ{t}{T}g(s,{Y}_{s})d{A}_s + \integ{t}{T} d{R}_s+\displaystyle{\sum_{t<s\leq
T}}h(s,Y_{s-} ,Y_s)\\
&\qquad\quad+\integ{t}{T}d{K}_{s}^+ -\integ{t}{T}d{K}_s^-
-\integ{t}{T}{Z}_{s}dB_{s}\,, t\leq T,
\\ (ii)&
\forall t\in[0,T[,\,\, L_t \leq Y_{t}\leq U_{t},
\\ (iii)&\integ{0}{T}({Y}_{t-}-{L}_{t-}) d{K}_{t}^+= \integ{0}{T}(
{U}_{t-}-{Y}_{t-}) d{K}_{t}^-=0,\,\, P-\mbox{a.s.}
\\ (iv)& {Y}\in {\cal D}, \quad {K}^+,
{K}^-\in {\cal K}, \quad {Z}\in {\cal L}^{2,d},
 \\ (v)&
d{K}^+\perp d{K}^-,
\end{array}
\right.
\end{equation} under the following assumptions :\\\\
\ni $(\bf{H.0})$ $A\in {\cal K}$ and $R\in {\cal K} - {\cal K}$ such
that: $0\leq A_T \leq 1$.

\ni $(\bf{H.1})$ There exist two processes $\eta \in
L^0(\Omega\times[0, T], \R_+))$ such that $\dint_0^T \eta_s ds\leq 1$
and $C\in \R_+ + {\cal K}$ such that:
\begin{enumerate}
 \item $\forall (s,\omega, y,z)\in [0,T]\times\Omega\times\R\times\R^d,\quad
-\eta_s(\omega)-{C_s(\omega)}|z|^2\leq f(s, \omega, y, z ) \leq
 0,$
 \item $\forall (t, \omega)\in [0, T]\times \Omega,$ the function $(y, z)\longmapsto f(t,
\omega, y, z)$ is continuous on $\R\times\R^d$.
\end{enumerate}
\ni $(\bf{H.2})$ For each $(s,\omega)\in [0,T]\times\Omega$,
  \begin{enumerate}
 \item $\forall y\in\R,\quad -1\leq g(s, \omega, y) \leq 0,$
 \item the function $ y\longmapsto g(s,
\omega, y)$ is continuous on $\R$.
 \end{enumerate}
\ni $(\bf{H.3})$ The function $h$ satisfies the following conditions
:
\begin{enumerate}
\item $P$-a.s.,
$\forall t\in]0, T]$,\,$\forall x\in\R$, the function $y\mapsto
{y}+h(t,x, {y}) $ is nondecreasing and continuous on $\R$.
\item $P$-a.s., $\forall t\in]0, T]$,\,$\forall y\in\R$, the function $x\mapsto
h(t,x, {y}) $ is continuous on $\R$.
\item There exists a nonnegative function $l : ]0, T]\times \Omega
\longrightarrow \R_+$ satisfying\, $\forall t\in ]0, T]$, $l_t$ is
${\cal F}_t$-measurable and $\displaystyle{\sum_{0< s\leq T}l_s}
<+\infty\,\,P-a.s.$ such that for every $t\in [0, T]$, $x\in\R$ and
$y\in\R$
$$
-l_t\leq h(t, x, {y})\leq 0.
$$
\end{enumerate}
\ni $(\bf{H.4})$ For each $(s,\omega)\in [0,T]\times\Omega$
$$
\begin{array}{ll}&
-1\leq L_s(\omega)\leq 0\leq U_s(\omega)\leq  1.
\end{array}
$$


Our main result Theorem \ref{thee1} is equivalent to the following
theorem.
\begin{theorem}\label{theo}
 Let assumptions $(\bf{H.0})$--$(\bf{H.4})$ hold. Then the GRBSDE
(\ref{eq221}) has a maximal solution.
\end{theorem}
\section{Proof of Theorem \ref{theo}}
For the proof of Theorem \ref{theo} we distinguish the following
four cases.
\subsection{Existence of solution for GBSDE (\ref{eq221}) : the "$f=g=h=0$" case}
\begin{theorem}\label{the11} Let assumptions $(\bf{H.0})$--$(\bf{H.4})$ hold. Assume moreover that $f=g=h=0$, then the GRBSDE
(\ref{eq221}) has a unique solution.
\end{theorem}
\ni In order to prove Theorem \ref{the11} we need some preliminary
results. To begin with, let $(\tau_n)_{n\geq 0}$ be the family of
stopping times defined by
\begin{equation}\label{tau1}
\tau_n =\inf\{s\geq 0 : |R|_s\geq n\}\wedge T.
\end{equation}
It easy seen that $P[\displaystyle{\cup_{n\geq 0}}(\tau_n =T)] =1$.
Indeed, let $\omega\in\displaystyle{\cap_{k\geq 0}} (\tau_k <T)$
then $\forall k,\, \tau_k <T$. Hence $\forall k,\,|R|_T\geq
k\Longleftrightarrow|R|_T =\infty$, which contradict the fact that
the total variation of $R$ is finite. \\ For all $(i,j)\in\N\times
\N$, let us set :
$$
\begin{array}{ll}
& \xi' = L_{T-}\vee (\Delta_T R)\wedge U_{T-}
\\ & L_t' = L_t  1_{\{t<T\}} + L_{T-} \, 1_{\{t=T\}}
\\ & U_t' = U_t  1_{\{t<T\}} + U_{T-} \, 1_{\{t=T\}}
\\ & \overline{\xi}^{i,j} = \xi'  +\dint_0^T 1_{\{s<\tau_i\}}
dR^{+}_{s}-\dint_0^T 1_{\{s<\tau_j\}} dR^{-}_{s}
\\ &
\overline{L}^{i,j}_t = L_t' +\dint_0^t 1_{\{s<\tau_i\}}
dR^{+}_{s}-\dint_0^t 1_{\{s<\tau_j\}} dR^{-}_{s}
\\ &
\overline{U}^{i,j}_t = U_t' +\dint_0^t 1_{\{s<\tau_i\}}
dR^{+}_{s}-\dint_0^t 1_{\{s<\tau_j\}} dR^{-}_{s}.
\end{array}
$$
Clearly we have $\overline{L}^{i,j}_t\leq \overline{\xi}^{i,j} \leq
\overline{U}^{i,j}_t$ and then by assumption $(\bf{H.4})$
$$
-(1+j)\leq \overline{L}^{i,j}_t\leq R^{i,j}_t:=\dint_0^t
1_{\{s<\tau_i\}} dR^{+}_{s}-\dint_0^t 1_{\{s<\tau_j\}}
dR^{-}_{s}\leq \overline{U}^{i,j}_t \leq 1+i.
$$
Consider the following BSDE with two reflecting barriers associated
with $(\overline{\xi}^{i,j}, \overline{L}^{i,j},
\overline{U}^{i,j})$
\begin{equation}
\label{eq7} \left\{
\begin{array}{ll}
(i) & 
 \overline{Y}_{t}^{i,j}=\overline{\xi}^{i,j} +\integ{t}{T}dK_{s}^{i,j+} -\integ{t}{T}dK_{s}^{i,j-}
-\integ{t}{T}Z_{s}^{i,j}dB_{s}\,, t\leq T,
\\ (ii)& \,\,
\forall  t\leq T,\,\, \overline{L}^{i,j}_t \leq \overline{Y}^{i,j}_t\leq \overline{U}^{i,j}_t,\\
(iii) & \integ{0}{T}(
\overline{Y}^{i,j}_{t-}-\overline{L}^{i,j}_{t-}) dK_{t}^{i,j+}=
\integ{0}{T}(
\overline{U}^{i,j}_{t-}-\overline{Y}^{i,j}_{t-}) dK_{t}^{i,j-}=0,\,\, \mbox{a.s.}\\
(iv)&\overline{Y}^{i,j}\in {\cal D}, \quad K^{i,j+}, K^{i,j-}\in
{\cal K}, \quad Z^{i,j}\in {\cal L}^{2,d},  \\ (v)& d{K}^{i,j+}\perp
d{K}^{i,j-}.
\end{array}
\right. \end{equation} It follows from Lepeltier and Xu \cite{Lexu}
(see Hamadène et \it{al.} \cite{HHO}) that Equation (\ref{eq7}) has
a unique solution. Moreover, for all $i$ and $j$
\begin{equation}\label{eqq4}
 \E\dint_0^T \mid Z^{i,j}_s\mid^2 ds+\E(K_T^{i,j\pm})^2 <+\infty.
\end{equation}
Set $Y^{i,j}_t =\overline{Y}^{i,j}_t -\dint_0^t 1_{\{s<\tau_i\}}
dR^{+}_{s}+\dint_0^t 1_{\{s<\tau_j\}} dR^{-}_{s}$ it follows then
that GBSDE (\ref{eq7}) can be written as follows:
\begin{equation}
\label{eq8} \left\{
\begin{array}{ll}
(i) & 
 Y_{t}^{i,j}=\xi'+\dint_t^T 1_{\{s<\tau_i\}}
dR^{+}_{s}-\dint_t^T 1_{\{s<\tau_j\}} dR^{-}_{s}
+\integ{t}{T}dK_{s}^{i,j+}\\ &
\qquad\qquad-\integ{t}{T}dK_{s}^{i,j-}
-\integ{t}{T}Z_{s}^{i,j}dB_{s},\,\, t\leq T,
\\ (ii)& \,\,
\forall  t\leq T,\,\, L_t' \leq Y^{i,j}_t\leq U_t',\\
(iii) & \integ{0}{T}( Y^{i,j}_{t-}-L_{t-}') dK_{t}^{i,j+}=
\integ{0}{T}(
U_{t-}'-Y^{i,j}_{t-}) dK_{t}^{i,j-}=0,\,\, \mbox{a.s.}\\
(iv)&Y^{i,j}_t\in {\cal D}, \quad K^{i,j+}, K^{i,j-}\in {\cal K},
\quad Z^{i,j}\in {\cal L}^{2,d},  \\ (v)& d{K}^{i,j+}\perp
d{K}^{i,j-}.
\end{array}
\right. \end{equation}

The following result follows easily from the Comparison theorem
(Theorem \ref{th111}. in Appendix).
\begin{proposition}\label{pro1} The solution $(Y^{i, j}, Z^{i, j}, K^{i, j+}, K^{i, j-})$ of RGBSDE (\ref{eq8}) satisfies the
following :
\begin{itemize}
\item[i)] Fix $j\in\N^*$, we get for all $i\geq 1$ and $t\leq T$
$$
-1\leq {L}_t'\leq Y_{t}^{i, j}\leq Y_{t}^{i+1, j}\leq {U}_{t}'\leq
1,\quad dK^{i+1,j+}\leq dK^{i,j+} \quad\mbox{ and }\quad
dK^{i,j-}\leq dK^{i+1, j-} .
$$
\item[ii)] Fix $i\in\N^*$, we get for all $j\geq
1$ and $t\leq T$
$$
{L}_t'\leq Y_{t}^{i, j+1}\leq Y_{t}^{i, j}\leq {U}_{t}',\quad
dK^{i,j+}\leq dK^{i,j+1+} \quad\mbox{ and }\quad dK^{j+1, i-}\leq
dK^{i,j-}.
$$
\end{itemize}
\end{proposition}
\ni Let us set\\
$\bullet$ $Y^j = \displaystyle\sup_{i}Y^{i, j}$,\,\, $Y^{j-}_{t} =
\displaystyle\sup_{i}Y^{i, j}_{t-}$.
\\
$\bullet$ $dK^{j-} = \displaystyle\sup_{i}dK^{i,j-}$ which is a
positive measure. \\
$\bullet$ $dK^{j+} = \displaystyle\inf_{i}dK^{i, j+}$ which is also
a positive measure since $K_T^{0, j+}<+\infty, \, P-a.s.$\\

The following result states the the existence of a process $Z^j$
such that process $(Y^j, Z^j, K^{j+}, K^{j-})$ is the unique
solution of some RBSDE.
\begin{proposition}\label{pro21} Suppose that assumptions of Theorem \ref{the11} hold. Then we have the following.
\begin{enumerate}
\item There exists a process $Z^j\in {\cal L}^{2,d}$ such that for
all $n\in\N$,
$$
\E\integ{0}{\tau_n}|Z_s^{i, j} - Z_s^{j}|^2ds \longrightarrow
0,\quad \mbox{as $i$ goes to infinity}.
$$
\item The process $(Y^j, Z^j, K^{j+}, K^{j-})$ is the unique solution
of the following BSDE with two reflecting barriers
\begin{equation}
\label{eqqq}
\left\{
\begin{array}{ll}
(i)&
 Y_{t}^{j}=\xi'+\dint_t^T 1_{\{s<T\}}
dR^{+}_{s}-\dint_t^T 1_{\{s<\tau_j\}} dR^{-}_{s}
+\integ{t}{T}dK_{s}^{j+}
\\ &\qquad-\integ{t}{T}dK_{s}^{j-} -\integ{t}{T}Z_{s}^{j}dB_{s}\,, t\leq T,
\\ (ii)& \,\,
\forall  t\leq T,\,\, L_t' \leq Y^{j}_t\leq U_t',\\
(iii) & \integ{0}{T}( Y^{j}_{t-}-L_{t-}') dK_{t}^{j+}= \integ{0}{T}(
U_{t-}'-Y^{j}_{t-}) dK_{t}^{j-}=0,\,\, \mbox{a.s.}\\
(iv)&Y^{j}_t\in {\cal D}, \quad K^{j+}, K^{j-}\in {\cal K}, \quad
Z^{j}\in {\cal L}^{2,d},  \\ (v)& d{K}^{j+}\perp d{K}^{j-}.
\end{array}
\right.
\end{equation}
\end{enumerate}
\end{proposition}\bop.
\ni \it{1.} Let $i, i'\in\N$ such that $i, i'\geq n$ and $t\in [0,
\tau_n[$ where $\tau_n$ is defined by (\ref{tau1}). Clearly we have
$$
\begin{array}{ll}
& (Y_{t}^{i, j} -Y_{t}^{i',j})^2 + \integ{t}{\tau_{n}}|Z_{s}^{i,j}-
Z_{s}^{i',j}|^2 ds\\ &= (Y_{\tau_{n-}}^{i, j}
-Y_{\tau_{n-}}^{i',j})^2 +2\integ{t}{\tau_{n-}}(Y_{s}^{i, j}
-Y_{s}^{i',j})((dK_{s}^{i,j+}-dK_{s}^{i',j+})-(dK_{s}^{i,j-}-dK_{s}^{i',j-}))
\\ &
-2\integ{t}{\tau_{n}}(Y_{s}^{i, j} -Y_{s}^{i',j})(Z_{s}^{i,j}-
Z_{s}^{i',j})dB_{s}- \sum_{t<s< \tau_n}(\Delta_s(Y^{i, j}
-Y^{i',j}))^2.
\end{array}
$$
By using a localization procedure it follows that
$$
\begin{array}{ll}
 \displaystyle{\lim_{i,i'\rightarrow +\infty}}\E\integ{0}{\tau_{n}}|Z_s^{i, j} - Z_s^{i',j}|^2ds
=0.
\end{array}
$$
Hence there exists $Z^j\in {\cal L}^{2,d}$ such that for every
$n\in\N$
$$
\begin{array}{ll}
\displaystyle{ \lim_{i\rightarrow
+\infty}}\E\integ{0}{\tau_{n}}|Z_s^{i, j} - Z_s^{j}|^2ds =0.
\end{array}
$$
2. According to Bulkholder-Davis-Gundy inequality there exists a
universal constant $c\geq 0$ such that
$$
\begin{array}{ll}
& \E\displaystyle{\sup_{t<\tau_n}}(Y_{t}^{i, j} -Y_{t}^{i',j})^2\leq
\E(Y_{\tau_{n-}}^{i, j} -Y_{\tau_{n-}}^{i',j})^2
+2c\E\bigg(\integ{0}{\tau_{n}}|Y_{s}^{i, j}
-Y_{s}^{i',j}|^2|Z_{s}^{i,j}- Z_{s}^{i',j}|^2ds\bigg)^{\frac12}.
\end{array}
$$
Henceforth
$$
\displaystyle{ \lim_{i\rightarrow
+\infty}}\E\displaystyle{\sup_{t<\tau_n}}(Y_{t}^{i, j} -Y_{t}^{j})^2
=0.
$$
Then $Y^j$ is \it{rcll} and $Y_t^{j-} = Y_{t-}^j, \forall t\in [0,
T],\,P-$a.s., since $P[\displaystyle{\cup_{n\geq 0}}(\tau_n =T)]
=1$. Moreover, since for every $(i,j)\in\N^2$,
$$
\E(K_T^{i,j+})^2\leq \E(K_T^{0,j+})^2 <+\infty,
$$
it follows then, from Fatou's lemma, that $\E(K_T^{j+})^2<+\infty$.\\
On the other hand, it follows from BSDE (\ref{eq8}) that for each
$i\geq n$ and $j\in\N$,
$$
\begin{array}{ll}
K_{\tau_{n-}}^{i,j-} &= Y_{\tau_{n-}}^{i,j}- Y_{0}^{i,j}+
K_{\tau_{n-}}^{i,j+}+\dint_0^{\tau_{n-}}1_{\{s<\tau_j\}}
dR^{+}_{s}-\dint_0^{\tau_{n-}} dR^{-}_{s}
 -\integ{0}{\tau_n}Z_{s}^{i,j}dB_{s}
 \\ & \leq 2+n+
K_{T}^{0,j+}
 -\integ{0}{\tau_n}Z_{s}^{i,j}dB_{s}.
 \end{array}
 $$
 Hence
$$
\begin{array}{ll}
\E K_{\tau_{n-}}^{i,j-} &\leq 2+n+\E K_{T}^{0,j+}.
 \end{array}
 $$
 By monotone convergence theorem we have
$$
\begin{array}{ll}
\E K_{\tau_{n-}}^{j-} &\leq 2+n+\E K_{T}^{0,j+}
 \end{array}
 $$
 Then $P-$a.s., for each $n\in\N$, $K_{\tau_{n-}}^{j-} <+\infty$.
 Therefore, since $P[\displaystyle{\cup_{n\geq 0}}(\tau_n =T)]
=1$ we get $K_{T-}^{j-}<+\infty$,\,\, $P-$a.s.
 Moreover
$$
\Delta_T K^{j-} = \sup_{i\in\N} \Delta_T K^{i,j-} = \sup_{i\in\N}
(U'_{T-}-(\xi'+\Delta_T R^{i,j}))^- =0.
$$
Henceforth $\forall j\in\N, K_{T}^{j-} <+\infty,\,\,P-$a.s. Similarly it follows also that $\Delta_T K^{j+} =0$. Therefore $K_{T}^{j\pm}<+\infty$,\,\, $P-$a.s.\\
Observe now that for every $n, j\in\N$, we have $P-$a.s.
$$
Y_{t}^{i,j}=Y_{\tau_{n-}}^{i,j}+\dint_t^{\tau_{n-}} (
1_{\{s<\tau_i\}}dR^{+}_{s}-1_{\{s<\tau_j\}} dR^{-}_{s})
+\integ{t}{\tau_{n-}}dK_{s}^{i,j+}-\integ{t}{\tau_{n-}}dK_{s}^{i,j-}
-\integ{t}{\tau_{n}}Z_{s}^{i,j}dB_{s}
$$
Taking the limit as $i$ goes to infinity we get for every $n,
j\in\N$, $P-$a.s.
$$
Y_{t}^{j}=Y_{\tau_{n-}}^{j}+\dint_t^{\tau_{n-}} (
dR^{+}_{s}-1_{\{s<\tau_j\}} dR^{-}_{s})
+\integ{t}{\tau_{n-}}dK_{s}^{j+}-\integ{t}{\tau_{n-}}dK_{s}^{j-}
-\integ{t}{\tau_{n}}Z_{s}^{j}dB_{s}
$$
Now, since $P[\displaystyle{\cup_{n\geq 0}}(\tau_n =T)] =1$, letting
$n$ goes to infinity we have for each $j\in\N$, $P-$a.s.
$$
Y_{t}^{j}=Y_{T-}^{j}+\dint_t^{T} (1_{\{s<T\}}
dR^{+}_{s}-1_{\{s<\tau_j\}} dR^{-}_{s})
+\integ{t}{T}dK_{s}^{j+}-\integ{t}{T}dK_{s}^{j-}
-\integ{t}{T}Z_{s}^{j}dB_{s}
$$
Therefore
$$
Y_{t}^{j} =\xi' +\dint_t^{T} (1_{\{s<T\}}
dR^{+}_{s}-1_{\{s<\tau_j\}} dR^{-}_{s})
+\integ{t}{T}dK_{s}^{j+}-\integ{t}{T}dK_{s}^{j-}
-\integ{t}{T}Z_{s}^{j}dB_{s},
$$
where we have used the fact that $Y^j_{T-} = \displaystyle{\sup_{i}}
Y^{i,j}_{T-} = \displaystyle{\sup_{i}} (L'_{T-}\vee
(\xi'+\Delta_TR^{i,j})\wedge U'_{T-}) =\xi'$, since $\Delta_TR^{i,j}
=0$. Then $(Y^{j}, Z^{j},K^{j+},K^{j-})$ satisfies $(i)$ of Equation
(\ref{eqqq}).\\

Let us now proof the minimality conditions. Clearly
$$
0\leq \integ{0}{T}(Y_{t-}^{i, j} -L'_{t-}) dK_{t}^{j+} \leq
\integ{0}{T}(Y_{t-}^{i, j} -L'_{t-}) dK_{t}^{i,j+} = 0.
$$
Hence
$$
\integ{0}{T}(Y_{t-}^{i, j} -L'_{t-}) dK_{t}^{j+} = 0.
$$
Applying Fatou's lemma we obtain
$$
0\leq\integ{0}{T}(Y_{t-}^{j} -L'_{t-}) dK_{t}^{j+} \leq \liminf_{i}
\integ{0}{T}(Y_{t-}^{i, j} -L'_{t-}) dK_{t}^{j+}= 0.
$$
Henceforth
$$
\integ{0}{T}(Y_{t-}^{j} -L'_{t-}) dK_{t}^{j+} = 0.
$$
Similarly we get also that
$$
\integ{0}{T}( U_{t-}'-Y^{j}_{t-}) dK_{t}^{j-}=0.
$$
Moreover, since $dK^{j+} = \displaystyle\inf_{i}dK^{i, j+}$,\,\,
$dK^{j-} = \displaystyle\sup_{i}dK^{i,j-}$ and $dK^{i, j+} \perp
dK^{i, j-}$ for each $i\in\N$ then $dK^{j+}\perp dK^{j-}$.
Proposition \ref{pro21} is then proved. \eop In view of passing to
the limit in Proposition \ref{pro1} (or in Theorem \ref{th111}. in
the Appendix), we get the following.
\begin{proposition}
For all $j\geq 1$, we obtain
\begin{enumerate}
\item $ Y^{j+1}\leq Y^j.$
\item $ dK_{t}^{j+} \leq dK_{t}^{j+1+}\,\,\,\, \mbox{and}\,\,\,\,\,
dK_{t}^{j+1-} \leq dK_{t}^{j-}.$
\end{enumerate}
\end{proposition}
 \ni Now let us set
\\
$\bullet$ $Y' = \displaystyle\inf_{j}Y^{j}$,\,\, $Y^{-}_t =
\displaystyle\inf_{j}Y^{j}_{t-}$.
\\
$\bullet$ $dK'^{-} = \displaystyle\inf_{j}dK^{j-}$ which is a
positive measure since $K_T^{0-}<+\infty, \, P-a.s.$ \\
$\bullet$ $dK'^{+} = \displaystyle\sup_{j}dK^{j+}$ which is also a
positive measure.\\
\begin{remark}
It should be noted that $Y'_T = Y_{T}^- = \xi'$.
\end{remark}
The following result states the convergence of the process $Z^j$ in
$L^2([0,\tau_n]\times \Omega)$.
\begin{proposition}\label{pro2} Suppose that assumptions of Theorem \ref{the11} hold. Then we have the
following:
\begin{enumerate}
\item There exists a process $Z\in {\cal L}^{2,d}$ such that, for
all $n\in\N$,
$$
\E\integ{0}{\tau_n}|Z_s^{j} - Z_s|^2ds \longrightarrow 0,\quad
\mbox{as $j$ goes to infinity}.
$$
\item The process $(Y', Z', K'^{+}, K'^{-})$ is the unique solution of the following GRBSDE with two
reflecting barriers
\begin{equation}
\label{eqq} \left\{
\begin{array}{ll}
(i) & 
 Y'_{t}={\xi'} +\dint_t^{T} 1_{\{s<T\}}
dR_s+\integ{t}{T}dK'^{+}_{s} -\integ{t}{T}dK'^{-}_{s}
-\integ{t}{T}Z_{s}dB_{s}\,, t\leq T,
\\ (ii)&
\forall  t\leq T,\,\, {L'}_t \leq Y'_{t}\leq {U'}_{t},\\
(iii) & \integ{0}{T}( Y'_{t-}-{L'}_{t-}) dK'^{+}_{t}= \integ{0}{T}(
{U'}_{t-}-Y'_{t-}) dK'^{-}_{t}=0,\,\, \mbox{a.s.}\\
(iv)&Y'\in {\cal D}, \quad K'^{+}, K'^{-}\in {\cal K}, \quad Z\in
{\cal L}^{2,d},
\\ (v)&
dK'^{+}\perp dK'^{-}.
\end{array}
\right. \end{equation}
\end{enumerate}
\end{proposition}
\bop. We just sketch the proof since the result follows by the same
way as previously. Let $n\in\N$ and  $j, j'\geq n$. By applying
It\^{o}'s formula to $(Y^j_t -Y^{j'}_t)^2$ on $[0,\tau_n[$ it
follows that there exists $Z\in {\cal L}^{2,d}$ such that
$$
\E\displaystyle{\sup_{t< \tau_n}}|Y'_t -Y_t^j|^2
+\E\integ{0}{\tau_n}|Z_s^{j} - Z_s|^2ds \longrightarrow 0,\quad
\mbox{as $j$ goes to infinity}.
$$
Hence $Y'$ is \it{rcll} and $Y^-_t =Y'_{t-},$\, $\forall t\in ]0,
T]\,P-$a.s. By the same way as previously we have also (\ref{eqq}).
\eop
\bop\,\,\bf{of Theorem \ref{the11}}. Uniqueness of solutions follows
easily. Let us focus on the existence. Let $Z$ be the process given
by Proposition \ref{pro2} and define
$$
\begin{array}{ll}
& Y_t = Y'_t  1_{\{t<T\}}
\\ & K_t^+ = K'^+_t + (L_{T-} -\Delta_T R)^+1_{\{t=T\}}
\\ & K_t^- = K'^-_t + (U_{T-}-\Delta_T R)^-1_{\{t=T\}}.
\end{array}
$$
Observe that for all $t\in [0, T]$
$$
Y_{t}=(\xi'-\Delta_T
R)1_{\{t<T\}}+\integ{t}{T}dR_s+\integ{t}{T}dK'^{+}_{s}
-\integ{t}{T}dK'^{-}_{s} -\integ{t}{T}Z_{s}dB_{s}.
$$
Since $Y'$ is \it{rcll} it follows also that $Y$ is \it{rcll}.
Clearly we have
$$
\begin{array}{ll}
\xi'-\Delta_T R &= (L_{T-}-\Delta_T R)\vee 0 \wedge (U_{T-}-\Delta_T
R)
\\ & = (L_{T-}-\Delta_T R)^+-(U_{T-}-\Delta_T R)^-,
\end{array}
$$
$$
\begin{array}{ll}
 K_T^+ -K_t^+ &= (K'^+_T + (L_{T-} -\Delta_T R)^+- K'^+_t)1_{\{t<T\}}
\\ & = K'^+_T -K'^+_t +(L_{T-} -\Delta_T R)^+1_{\{t<T\}},
\end{array}
$$
$$
\begin{array}{ll}
K_{T}^- -K_t^- = K'^-_T -K'^-_t +(U_{T-} -\Delta_T R)^+1_{\{t<T\}}.
\end{array}
$$
Consequently  $$ \left\{
\begin{array}{ll}
 & 
 Y_{t}=\dint_t^{T} dR_{s} +\integ{t}{T}dK_{s}^{+}
-\integ{t}{T}dK_{s}^{-} -\integ{t}{T}Z_{s}dB_{s}\,, t\leq T,
\\ &
\forall  t\leq T,\,\, {L}_t \leq Y_{t}\leq {U}_{t}.
\end{array}
\right.
$$
Moreover it follows that 
$$
\begin{array}{ll}
\integ{0}{T}( Y_{t-}-{L}_{t-}) dK_{t}^{+} &= ((L_{T-} \vee (\Delta_T
R)\wedge U_{T-}) -L_{T-})(L_{T-}-\Delta_T R)^+
\\ & = ((0 \vee (-L_{T-}+\Delta_T R)\wedge (U_{T-} -L_{T-}))(L_{T-}-\Delta_T
R)^+
\\ & = 0.
\end{array}
$$
Similarly we get also that
$$
\integ{0}{T}( U_{t-}-{Y}_{t-}) dK_{t}^{-}=0.
$$
Moreover, since $L_{T-}\leq U_{T-}$, we have $ (L_{T-}-\Delta_T
R)^+(U_{T-}-\Delta_T R)^- =0$,  then $dK^{+}\perp dK^{-}$. The proof
of Theorem \ref{the11} is finished.\eop

\subsection{Existence of solution for GBSDE (\ref{eq221}): the "$f$, $g$ are
Lipschitz and $h=0$" case}
\begin{theorem}\label{the1} Let assumptions $(\bf{H.0})$--$(\bf{H.4})$ hold. Assume moreover that
$h=0$ and $f$ and $g$ are $a-$Lipschitz, then the GRBSDE
(\ref{eq221}) has a unique solution.
\end{theorem}
\bop. The existence proof is based on the Picard's approximation
scheme. Let $(Y^0, Z^0, K^{+,0}, K^{,0}) =(0, 0, 0, 0)$ and define
$(Y^{n+1}, Z^{n+1}, K^{n+1+}, K^{n+1-})$ as the solution (which
exists according to the previous subsection) of the following GRBSDE
\begin{equation}
\label{eq91} \left\{
\begin{array}{ll}
(i) & 
 {Y}_{t}^{n+1}=\integ{t}{T}{f}(s,{Y}_{s}^n,{Z}_{s}^n)ds +
 \integ{t}{T}g(s,{Y}_{s}^n)d{A}_s + \integ{t}{T} d{R}_s\\
&\qquad\quad+\integ{t}{T}d{K}_{s}^{n+1+} -\integ{t}{T}d{K}_s^{n+1-}
-\integ{t}{T}{Z}_{s}^{n+1}dB_{s}\,, t\leq T,\\ (ii)& \forall
t\in[0,T],\,\, L_{t} \leq Y_{t}^{n+1} \leq U_{t},
\\ (iii)&\integ{0}{T}({Y}_{t-}^{n+1}-{L}_{t-}) d{K}_{t}^{n+1+}= \integ{0}{T}(
{U}_{t-}-{Y}_{t-}^{n+1}) d{K}_{t}^{n+1-}=0,\,\, P-\mbox{a.s.}
\\ (iv)& {Y}^{n+1}\in {\cal D} \quad {K}^{n+1+},
{K}^{n+1-}\in {\cal K} \quad {Z}^{n+1}\in {\cal L}^{2,d},
 \\ (v)&
d{K}^{n+1+}\perp d{K}^{n+1-},
\end{array}
\right.
\end{equation}
Let $n, m\in\N$. By applying It\^{o}'s formula to
$(Y_t^{n+1}-Y_t^{m+1})^2e^{\alpha(t+A_t)},\,\, \alpha\geq
8a^2(T+1)(1+4c^2)$, where $c$ is a universal constant, coming from
Bulkholder-Davis-Gundy inequality, and using standard calculations
for RBSDE one can prove that there exists a process $ Y\in {\cal D}$
and $Z\in {\cal L}^{2,d}$ such that
$$
\lim_{n\rightarrow +\infty}\E\sup_{t\leq T}(Y_{t}^{n} -Y_{t})^2
=0,\quad \lim_{n\rightarrow +\infty}\E\integ{0}{T}|Z_s^{n} -
Z_s|^2ds = 0.
$$
Let $(\overline{Y}, \overline{Z}, K^+, K^-)$ be the solution, which
exists according to the previous subsection, of the following GRBSDE
\begin{equation}
\left\{
\begin{array}{ll}
(i) & \overline{Y}_{t}=\integ{t}{T}{f}(s, Y_{s},Z_{s})ds +
 \integ{t}{T}g(s,{Y}_{s})d{A}_s + \integ{t}{T} d{R}_s\\
&\qquad\quad+\integ{t}{T}d{K}_{s}^{+} -\integ{t}{T}d{K}_s^{-}
-\integ{t}{T}\overline{Z}_{s}dB_{s}\,, t\leq T,
\\ (ii)&
\forall  t\leq T,\,\, {L}_t \leq \overline{Y}_{t}\leq {U}_{t},\\
(iii) &  \integ{0}{T}(\overline{Y}_{t-}-{L}_{t-}) dK_{t}^{+}=
\integ{0}{T}( {U}_{t-}-\overline{Y}_{t-}) d{K}_{t}^{-}=0,\,\,
\mbox{a.s.}
\\ (iv)& \overline{Y}\in {\cal D}, \quad {K}^{+},
{K}^{-}\in {\cal K},\quad \overline{Z}\in {\cal L}^{2,d},
 \\ (v)&
d{K}^{+}\perp d{K}^{-}.
\end{array}
\right.
\end{equation}
It is not difficult to prove that there exists a constant $C>0$ such
that
$$
\E\sup_{t\leq T}(Y_{t}^{n+1} -\overline{Y}_{t})^2
+\E\integ{0}{T}|Z_s^{n+1} - \overline{Z}_s|^2ds\leq C (\E\sup_{t\leq
T}(Y_{t}^{n} -Y_t)^2+\E\integ{0}{T}|Z_s^{n} - Z_s|^2ds).
$$
Hence
$$
\lim_{n\rightarrow +\infty}[\E\sup_{t\leq T}(Y_{t}^{n}
-\overline{Y}_{t})^2 + \E\integ{0}{T}|Z_s^{n+1} -
\overline{Z}_s|^2ds]=0.
$$
It follows that
$$
\E\sup_{s\leq T}\mid Y_s -\overline{Y}_s \mid^2 =
0,\,\,\mbox{and}\,\, \E\integ{0}{T}|Z_s - \overline{Z}_s|^2ds =0.
$$
Therefore $Y = \overline{Y}$ and $Z=\overline{Z}$. The proof of
existence is then
 finished. The uniqueness of solutions follows easily by using standard
arguments. \eop 
\subsection{Existence of solution for GBSDE (\ref{eq221}) : the "$f$, $g$ are Lipschitz and there exists $S_0 =0\leq S_1\leq...\leq S_{p+1}=T$ such that
$\forall x,y\in\R$, $\forall t\notin \{S_1,...,S_{p+1}\}$ $h(t,
\omega,
x, y)=0$" case.} 
\begin{theorem}\label{the1} Let assumptions $(\bf{H.0})$--$(\bf{H.4})$ hold. Assume moreover
that there exists a finite family of stopping times $S_0 =0\leq
S_1\leq...\leq S_{p+1}=T$ such that for each $x,y\in\R$ and $t\notin
\{S_1,...,S_{p+1}\}$ $h(t, \omega, x, y)=0$ and $f$ and $g$ are
$a-$Lipschitz, then the GRBSDE (\ref{eq221}) has a maximal solution.
\end{theorem}
\bop. For every $y\in\R$ and $i\in\{1, 2,...,p+1\}$, set
$$
\begin{array}{ll}
&\widehat{h}(t, \omega, y) = \max\{ x\in[L_{t-}, U_{t-}]\,\,
: \,\, x= L_{t-}\vee[ y+ h(t, x, y)+\Delta_t R]\wedge U_{t-} \} \\
\\ & L_t^i = L_t  1_{\{t<S_i\}} + L_{S_{i-}}  1_{\{t\geq S_i\}}
\\ & U_t^i = U_t  1_{\{t<S_i\}} + U_{S_{i-}}   1_{\{t\geq S_i\}}
\\ & \xi^{p+1} = \widehat{h}(T, \omega, 0) = \widehat{h}(S_{p+1}, \omega,0).
\end{array}
$$
Let $(Y^{p+1}, Z^{p+1}, K^{p+1+}, K^{p+1-})$ be the unique solution
(which exists according to the previous subsection) of following
GRBSDE
\begin{equation}\label{eqqq1}
\left\{
\begin{array}{ll}
(i) & Y_{t}^{p+1}={\xi}^{p+1} +\integ{t}{T}1_{\{s<S_{p+1}\}}{f}(s,
Y_{s}^{p+1},Z_{s}^{p+1})ds +
 \integ{t}{T}1_{\{s<S_{p+1}\}}g(s,{Y}_{s}^{p+1})d{A}_s\\ & \qquad\quad+ \integ{t}{T} 1_{\{s<S_{p+1}\}} d{R}_s
+\integ{t}{T}d{K}_{s}^{p+1+} -\integ{t}{T}d{K}_s^{p+1-}
-\integ{t}{T}Z_{s}^{p+1}dB_{s}\,, t\leq T,
\\ (ii)&
\forall  t\leq T,\,\, {L}_t^{p+1} \leq Y_{t}^{p+1}\leq {U}_{t}^{p+1},\\
(iii) &  \integ{0}{T}(Y_{t-}^{p+1}-{L}_{t-}^{p+1}) dK_{t}^{p+1+}=
\integ{0}{T}( {U}_{t-}^{p+1}-Y_{t-}^{p+1}) d{K}_{t}^{p+1-}=0,\,\,
\mbox{a.s.}
\\ (iv)& Y^{p+1}\in {\cal D}, \quad {K}^{p+1+},
{K}^{p+1-}\in {\cal K},\quad Z^{p+1}\in {\cal L}^{2,d},
 \\ (v)&
d{K}^{p+1+}\perp d{K}^{p+1-}.
\end{array}
\right.
\end{equation}
We should remark here, by Lemma \ref{lem1}, that
$$ Y_{S_{p+1}^-} =
L_{S_{p+1}^-} \vee Y_{S_{p+1}}\wedge U_{S_{p+1}^-} = Y_{S_{p+1}}=
\xi^{p+1} = \widehat{h}(S_{p+1}, \omega,0).
$$
Therefore $\Delta_{S_{p+1}}Y^{p+1} =0. $
\\
Now we want to construct a solution to our GRBSDE by induction and
concatenation. For that reason, suppose there exists a solution to
our GRBSDE $(Y^{i+1}, Z^{i+1}, K^{i+1+}, K^{i+1-})$ on $[0,T]$, for
$i\in \{1,...,p\}$. Let  $\xi^{i} = \widehat{h}(S_i, \omega,
Y_{S_i}^{i+1})$ and define $(Y^{i}, Z^{i}, K^{i+}, K^{i-})$ as the
unique solution of the following GRBSDE
\begin{equation}\label{eqqq2}
\left\{
\begin{array}{ll}
(i) & Y_{t}^{i}=\xi^i+\integ{t}{T}1_{\{s<S_{i}\}}{f}(s,
Y_{s}^{i},Z_{s}^{i})ds +
 \integ{t}{T}1_{\{s<S_{i}\}} g(s,{Y}_{s}^{i})d{A}_s + \integ{t}{T} 1_{\{s<S_{i}\}} d{R}_s\\
&\qquad\quad+\integ{t}{T}d{K}_{s}^{i+} -\integ{t}{T}d{K}_s^{i-}
-\integ{t}{T}Z_{s}^{i}dB_{s}\,, t\leq T,
\\ (ii)&
\forall  t\leq T,\,\, {L}_t^{i} \leq Y_{t}^{i}\leq {U}_{t}^{i},\\
(iii) &  \integ{0}{T}(Y_{t-}^{i}-{L}_{t-}^{i}) dK_{t}^{i+}=
\integ{0}{T}( {U}_{t-}^{i}-Y_{t-}^{i}) d{K}_{t}^{i-}=0,\,\,
\mbox{a.s.}
\\ (iv)& Y^i\in {\cal D}, \quad {K}^{i+},
{K}^{i-}\in {\cal K},\quad Z^{i}\in {\cal L}^{2,d},
 \\ (v)&
d{K}^{i+}\perp d{K}^{i-}.
\end{array}
\right.
\end{equation}
We have $\Delta_{S_{i}}Y^{i} =0$ and $Y^i_s= Y^i_{S_i}=\xi^i$,
$\forall s\in [S_i, T]$. By setting
$$
\begin{array}{ll}
& Y_t = \displaystyle{\sum_{i=0}^{p}} Y_t^{i+1}1_{[S_i, S_{i+1}[}(t)
\\ & Z_t = \displaystyle{\sum_{i=0}^{p}} Z_t^{i+1}1_{[S_i, S_{i+1}[}(t)
\\ & K_t^{\pm c} = \displaystyle{\sum_{i=0}^{p}} \dint_0^t 1_{[S_i,
S_{i+1}[}(s)dK_s^{i+1\pm c},
\end{array}
$$
it follows that
\begin{equation}
\label{eqqq3} \left\{
\begin{array}{ll}
(i) & 
 Y_{t}=\integ{t}{T}f(s,Y_{s},Z_{s})ds+\dint_t^Tg(s,
Y_s)dA_s+\integ{t}{T} d{R}_s^c -\sum_{t<s\leq T}\Delta Y_s\\
&\qquad\quad+\integ{t}{T}dK_{s}^{+c} -\integ{t}{T}dK_{s}^{-c}
-\integ{t}{T}Z_{s}dB_{s}\,, t\leq T,
\\ (ii)& \integ{0}{T}( Y_{t-}-L_{t-})
dK_{t}^{+c}= \integ{0}{T}( U_{t-}-Y_{t-}) dK_{t}^{-c}=0,\,\, \mbox{a.s.}, \\
(iii)& Y\in {\cal D}, \quad K^{+c}, K^{-c}\in {\cal K}^c, \quad Z\in
{\cal L}^{2,d},  \\ (v)& dK^{+c}\perp  dK^{-c}.
\end{array}
\right. \end{equation} Let $t\in ]0, T]$. We discus the following two cases. \\
1. If $t\notin \{S_1,...,S_{p+1}\}$, then there exists $i$ such that
$t\in ]S_i, S_{i+1}[$ and then, since $h(t,x,y) =0$, for every
$x,y$, we have
$$
\begin{array}{ll}
Y_{t-}&=Y_{t-}^{i+1} = L_{t-}^{i+1} \vee [Y_t^{i+1} + \Delta_t
R]\wedge U_{t-}
\\ & = \max\{ x\in[L_{t-}, U_{t-}]\,\,
: \,\, x= L_{t-}\vee[ Y_t+ h(t,x,Y_t)+\Delta_t R]\wedge U_{t-} \}.
\end{array}
$$
2. If there exits $i$ such that $t= S_i$, then
$$
\begin{array}{ll}
Y_{{S_i}^-}&=Y_{{S_i}^-}^{i} = Y^i_{S_i} =\xi^i
\\ & = \max\{ x\in[L_{{S_i}^-}, U_{{S_i}^-}]\,\,
: \,\, x= L_{{S_i}^-}\vee[ Y_{S_i}+h(S_i, x, Y_{S_i})+ \Delta_{S_i}
R]\wedge U_{{S_i}^-} \}.
\end{array}
$$
Hence for each $t\in ]0, T]$,
$$
\begin{array}{ll}
Y_{t-}= \max\{ x\in[L_{t-}, U_{t-}]\,\, : \,\, x= L_{t-}\vee[ Y_t+
h(t,x,Y_t)+\Delta_t R]\wedge U_{t-} \}.
\end{array}
$$
Consequently, by Lemma \ref{lem1}, $(Y, Z, K^+, K^-)$ is the unique
solution of the following GRBSDE
\begin{equation}
\label{eqqq3} \left\{
\begin{array}{ll}
(i) & 
 Y_{t}=\xi
+\integ{t}{T}f(s,Y_{s},Z_{s})ds+\dint_t^Tg(s,
Y_s)dA_s+\integ{t}{T}d{R}_s +\sum_{t<s\leq T}h(s, Y_{s-}, Y_s)\\
&\qquad\quad+\integ{t}{T}dK_{s}^{+} -\integ{t}{T}dK_{s}^{-}
-\integ{t}{T}Z_{s}dB_{s}\,, t\leq T,
\\ (ii)&
\forall t\leq T,\,\, Y_{t-} = \max\{ x\in[L_{t-}, U_{t-}]\,\, : \,\,
x= L_{t-}\vee[
Y_t+ h(t,x,Y_t)+\Delta_t R]\wedge U_{t-} \},\\
\\ (iii)& \integ{0}{T}( Y_{t-}-L_{t-})
dK_{t}^{+}= \integ{0}{T}( U_{t-}-Y_{t-}) dK_{t}^{-}=0,\,\, \mbox{a.s.}, \\
(iv)& Y\in {\cal D}, \quad K^{+}, K^{-}\in {\cal K}, \quad Z\in
{\cal L}^{2,d},  \\ (v)& dK^{+}\perp  dK^{-}.
\end{array}
\right. \end{equation} Henceforth $(Y, Z, K^+, K^-)$ is the maximal
solution of the GRBSDE (\ref{eq221}).\eop

\subsection{Existence of solution for GBSDE (\ref{eq221}) : the "general" case}\label{sec4.4}
\begin{theorem}\label{the122} Let assumptions $(\bf{H.0})$--$(\bf{H.4})$ hold.
Then the GRBSDE (\ref{eq221}) has a maximal solution.
\end{theorem}
Since the integrability conditions on parameters are weaker, the
proof of Theorem \ref{the122} is based on regularization by
sup-convolution techniques and a truncation procedure by means of a
family of stopping times. The final step consists in justifying the
passage to the limit and identifying the limit as the solution of
our GRBSDE (\ref{eq221}).
\subsubsection{Approximations}
It is not difficult to prove the following lemma which gives an
approximation of continuous functions by Lipschitz functions.
\begin{lemma}\label{lem0}Let $(T_i)_{i\geq 1}$ be a sequence of stopping times such that $[|T_i|]\cap[|T_j|] =
\emptyset,\,\, \forall i\neq j$ and $\bigcup_{i\geq 1}[|T_i|] =\{(t,
\omega)\in ]0,T[\times\Omega : l_t
>0\}$. For every $n\in\N$ define the
functions $f_n$, $g_n$ and $h_n$ by :
$$
\begin{array}{lll}
 &f_n(t,y,z) = \displaystyle\sup_{p\in \R, q\in \R^d} \{ {f}(t, p, q)-n|p-y| - n|q-z| \}1_{\{n\geq 1\}},
\\ &
g_n(t,y) = \displaystyle\sup_{p\in \R} \{{g}(t, p)-n|p-y|
\}1_{\{n\geq 1\}}
\\ &
h_n(t,x,y) =  {h}(t, x, y)1_{\{t\in\{ T_1,...,T_n\}\}}1_{\{n\geq
1\}}.
\end{array}
$$
 Assume that assumptions $(\bf{H.0})$--$(\bf{H.4})$ hold. Then we
have the following :
\begin{enumerate}
\item For all $(t,\omega, y,z,n)\in [0,T]\times\Omega\times\R\times\R^d\times\N,\,\,$
$$f_0(t,y,z)=0\geq f_n(t,y,z)
\geq f_{n+1}(t,y,z)\geq {f}(t, y, z) \geq -{\eta}_t -{C_t}|z|^2.$$
\item For all $(t,\omega, y,n)\in [0,T]\times\Omega\times\R\times\N,\,\,$
$$g_0(t,y)=0\geq g_n(t,y) \geq g_{n+1}(t,y)\geq {g}(t, y) \geq -1.$$
\item For all $(t,\omega,y,x, n)\in [0,T]\times\Omega\times\R\times\R\times\N,$
$$
 h_{0}(t,x, y) = 0\geq h_{n}(t,x,
y)\geq h_{n+1}(t,x, y)\geq h(t,x, y)\geq -l_t.
$$
\item $f_n$ is uniformly $n$-Lipschitz with respect to $(y,z)$.
\item $g_n$ is uniformly $n$-Lipschitz with respect to $y$.
\item For all $(t,\omega)\in [0,T]\times\Omega$, $ (f_n(t,y,z)) _{n\geq 0}$ converges to ${f}(t, y, z)$
as $n$ goes to $+\infty$ uniformly on every compact of $\R\times
\R^d.$
\item For all $(t,\omega)\in [0,T]\times\Omega$, $ (g_n(t,y)) _{n\geq 0}$ converges
to ${g}(t, y)$ as $n$ goes to $+\infty$ uniformly on every compact
of $\R$.
\item For all $(t,\omega)\in [0,T]\times\Omega\times\R\times\R$, $ (h_n(t,x,y)) _{n\geq 0}$ converges
to ${h}(t,x, y)$ as $n$ goes to $+\infty$ uniformly on every compact
of $\R\times\R$.
\end{enumerate}
\end{lemma}
Let us define for $i\in\{1,...,n-1\}$
$$
\begin{array}{ll}
& S_0^n =0
\\ &S_1^n = \min\{T_1,..., T_n\}
\\ & S_{i+1}^{n} = \min(\{T_1,..., T_n\}\backslash \{S_1^n,...,S_i^n\})
\\ & S_{n+1}^n = T.
\end{array}
$$
We note here that for $i\in\{1,...,n\}$, $S_{i}^{n}$ is a stopping
time. Indeed, it is clear that $S_1^n$ is a stopping time. Suppose
that for $i\in\{1,...,n-1\}$, $S_{i}^{n}$ is a stopping time and
prove that $S_{i+1}^{n}$ is a stopping time which is evident since
for every $t\in [0, T]$,
$$
\{S_{i+1}^{n}\leq t\} = \displaystyle{\bigcup_{j=1}^{n}}
\{S_{i}^{n}<T_j\leq t\}\in {\cal F}_t.
$$
Observe that
$$
\{T_1,..., T_n\}= \{S_1^n,...,S_n^n\}\quad \mbox{and}\quad 0=
S_0^n<S_1^n<...<S_n^n\leq S_{n+1}^n =T.
$$
 Henceforth
$$
h_n(t,x,y) =0,\,\, \forall t\notin\{ S_1^n,...,
S_{n+1}^n\},\,\forall (x,y)\in\R\times\R.
$$


Therefore, according to the previous subsection, for each $n\in\N$
there exists a unique solution $(Y^n, Z^n, K^{n+}, K^{n-})$ to the
following GRBSDE
\begin{equation}
\label{eq9} \left\{
\begin{array}{ll}
(i) & 
 {Y}_{t}^n={\xi} +\integ{t}{T}{f}_n(s,{Y}_{s}^n,{Z}_{s}^n)ds +
 \integ{t}{T}g_n(s,{Y}_{s}^n)d{A}_s + \integ{t}{T} d{R}_s\\
&\qquad\quad+\displaystyle{\sum_{t<s\leq T}}h_n(s,Y_{s-}^n,
Y_s^n)+\integ{t}{T}d{K}_{s}^{n+} -\integ{t}{T}d{K}_s^{n-}
-\integ{t}{T}{Z}_{s}^ndB_{s}\,, t\leq T,
\\ (ii)&
\forall t\leq T,\,\, Y_{t-}^n = \max\{ x\in[L_{t-}, U_{t-}]\,\, :
\,\, x= L_{t-}\vee[ Y_t^n+ h_n(t,x,Y_t^n)+\Delta_t R]\wedge U_{t-}
\},
\\ (iii)&\integ{0}{T}({Y}_{t-}^n-{L}_{t-}) d{K}_{t}^{n+}= \integ{0}{T}(
{U}_{t-}-{Y}_{t-}) d{K}_{t}^{n-}=0,\,\, \mbox{a.s.}
\\ (iv)& {Y}^n\in {\cal D}, \quad {K}^{+n},
{K}^{-n}\in {\cal K}^c, \quad {Z}^n\in {\cal L}^{2,d},
 \\ (v)&
d{K}^{n+}\perp d{K}^{n-}.
\end{array}
\right.
\end{equation}
\subsubsection{Convergence of the approximating scheme}\label{sub4.4.2} By using Comparison Theorem (Theorem \ref{th111}.) it is not
difficult to prove the following proposition.
\begin{proposition} For all
$n\geq 0$, we obtain

\begin{enumerate}
\item $L_t\leq Y^{n+1}_t\leq Y^n_t\leq U_t,\,\, \forall t\in [0,T],
P-a.s.$ \item $dK_{t}^{n+} \leq dK_{t}^{n+1+}\,\,\,\,
\mbox{and}\,\,\,\,\, dK_{t}^{n+1-} \leq
dK_{t}^{n-}\,\,\mbox{on}\,\,[0, T].$
\end{enumerate}
\end{proposition}
\ni Set
\\
$\bullet$ $Y_t = \displaystyle\inf_{n}Y^{n}_t, \,\, Y_t^{-} =
\displaystyle\inf_{n}Y^{n}_{t-}$.
\\
$\bullet$ $dK^{-}_t = \displaystyle\inf_{n}dK^{n-}_t$ which is a
positive measure on $[0,T]$ since $K_T^{0-}<+\infty\,\,P-$a.s. \\
$\bullet$ $dK^{+}_t = \displaystyle\sup_{j}dK^{n+}_t$ which is a
positive measure on $[0,T]$.\\

\ni Let $(\tau_j)_{j\geq 1}$ be the family of stopping times defined
by
\begin{equation}\label{tau2}
\tau_j =\inf\{s\geq 0 : \displaystyle{\sum_{r\leq s}}l_r+C_s
+|R|_s\geq j\}\wedge T.
\end{equation}
It should be pointed out that $P[\displaystyle{\cup_{j\geq
1}}(\tau_j =T)]
=1$.\\

The following result states the convergence of the process $Z^n$ in
$L^2([0,\tau_j]\times \Omega)$.
\begin{proposition}\label{th1} Assume that assumptions $(\bf{H.0})$--$(\bf{H.4})$
hold. Then there exists a process $Z\in {\cal L}^{2,d}$ such that,
for all $j$,
$$
\E\integ{0}{\tau_j}|Z_s^{n} - Z_s|^2ds \longrightarrow 0,\quad
\mbox{as $n$ goes to infinity}.
$$
\end{proposition}
\bop. For $s\in\R$ and $j\in \N$, let us set $\psi(s) = \dfrac{e^{4j
s}-1}{4j} -s$. We mention that $\psi$ satisfies the following for
all  $s\in \R$,
\begin{equation} \label{eqq7}
\begin{array}{lll}
 & \psi'(s) = e^{4j s} -1 ,\,\,\,  \psi''(s) = 4j e^{4j s} = 4j
 \psi'(s)+4j.
\end{array}
\end{equation}
Let $n,m\in \N$ such that $m\geq n$. Applying It\^{o}'s formula to
$\psi(Y^n -Y^m)$, we get for $t< \tau_j$,
$$
\begin{array}{ll}
&\psi(Y_{t}^n -Y_{t}^{m})
\\&=\psi(Y_{{\tau_j}^-}^n-Y_{{\tau_j}^-}^m)+
\integ{t}{{\tau_j}}(f_n(s,Y_{s}^n,Z_{s}^n)
-f_m(s,Y_{s}^{m},Z_{s}^{m}))\psi'(Y_{s}^n -Y_{s}^{m})ds\\
& +\integ{t}{{\tau_j}} (g_n(s,Y_{s}^{n})-g_m(s,Y_{s}^{m}))
\psi'(Y_{s}^n -Y_{s}^{m})dA_s
\\
&+\displaystyle{\sum_{t<s< {\tau_j}}}\psi'(Y_{s-}^n
-Y_{s-}^{m})(h_n(s, Y_{s-}^{n},Y_s^n)-h_m(s,Y_{s-}^{m}, Y_s^m))
\\ &
+\integ{t}{{\tau_j}^-}\psi'(Y_{s-}^n -Y_{s-}^{m})
d(K_{s}^{n+}-K_{s}^{m+})
-\integ{t}{{\tau_j}^-}\psi'(Y_{s-}^{n} -Y_{s-}^{m}) d(K_{s}^{n-}-K_{s}^{m-}) \\
&-\integ{t}{{\tau_j}}\psi'(Y_{s}^n -Y_{s}^{m})(Z_{s}^n-
Z_{s}^{m})dB_{s} -\dfrac12\integ{t}{{\tau_j}}\psi''(Y_{s}^n
-Y_{s}^{m})|Z_{s}^n- Z_{s}^{m}|^2 ds\\
& -\displaystyle{\sum_{t<s<{\tau_j}}}\bigg[\psi(Y_{s}^n
-Y_{s}^{m})-\psi(Y_{s-}^n -Y_{s-}^{m})-\psi'(Y_{s-}^n -Y_{s-
}^{m})\Delta(Y_{s}^n -Y_{s}^{m}))\bigg].
\end{array}
$$
Since $\psi'(0) =0$, we have $$ \integ{t}{{\tau_j}^-}\psi'(Y_{s-}^n
-Y^m_{s-}) d(K_{s}^{n+}-K_{s}^{m+}) =
-\integ{t}{{\tau_j}^-}\psi'(Y_{s-}^n -{L}_{s-}) dK_{s}^{m+}\leq 0,$$
and
 $$
-\integ{t}{{\tau_j}^-}\psi'(Y_{s-}^n -Y^m_{s-})
d(K_{s}^{n-}-K_{s}^{m-}) =
-\integ{t}{{\tau_j}^-}\psi'({U}_{s-}-Y_{s-}^m) dK_{s}^{n-}\leq 0.$$
Then
\begin{equation}\label{eqqq41}
\begin{array}{ll}
&\psi(Y_{t}^n -Y_{t}^{m})
\\&=\psi(Y_{{\tau_j}^-}^n-Y_{{\tau_j}^-}^m) -\integ{t}{{\tau_j}^-} dR^{n, m,j}_s +
\integ{t}{{\tau_j}}(\eta_s +j|Z_s^m|^2)\psi'(Y_{s}^n -Y_{s}^{m})ds\\
& +\integ{t}{{\tau_j}} \psi'(Y_{s}^n -Y_{s}^{m})dA_s
+\displaystyle{\sum_{t<s< {\tau_j}}}\psi'(Y_{s-}^n -Y_{s-}^{m})\,
l_s -\integ{t}{{\tau_j}}\psi'(Y_{s}^n -Y_{s}^{m})(Z_{s}^n-
Z_{s}^{m})dB_{s}\\
& -2j\integ{t}{{\tau_j}}\psi'(Y_{s}^n -Y_{s}^{m})|Z_{s}^n-
Z_{s}^{m}|^2 ds -2j\integ{t}{{\tau_j}}|Z_{s}^n- Z_{s}^{m}|^2 ds.
\end{array}
\end{equation}
where $dR^{n, m,j}$ is a positive measure depending on $n, m$ and
$j$. By taking $n=0$ in Equation (\ref{eqqq41}) and using a
localization procedure, we get for all $j\in\N$
\begin{equation}\label{eqqq5}
\begin{array}{ll}
&2j\E\integ{0}{{\tau_j}}|Z_{s}^0- Z_{s}^{m}|^2 ds
\\&\leq \psi(2)  +
\psi'(2)\E\integ{0}{T}\eta_sds
+2j\psi'(2)\E\integ{0}{\tau_j}|Z_s^0|^2ds+\psi'(2)\E A_T
+\psi'(2)\E\displaystyle{\sum_{0<s< {\tau_j}}}l_s.
\end{array}
\end{equation}
Hence there exists a positive constant  $c_j$ depending only on $j$
such that
 $$
\begin{array}{ll}
& \E\integ{0}{\tau_j}|Z_{s}^0- Z_{s}^{m}|^2 ds \leq c_{j}.
\end{array}
$$ 
It follows from subsection 4.1 that for all $j\geq 1$
$\E\integ{0}{\tau_j}|Z_{s}^0|^2 ds <+\infty $, then for all $j
\in\N^*$
$$
\begin{array}{ll}
& \displaystyle{\sup_{m\in\N}}\E\integ{0}{\tau_j}|Z_{s}^{m}|^2 ds
<+\infty.
\end{array}
$$
Henceforth, there exist a subsequence $\bigg(m_k^j\bigg)_{k}$ of $m$
and a process $\widehat{Z}^j\in L^2(\Omega\times [0,T]; \R^d))$ such
that
\\
$Z_s^{m_k^j}1_{\{s\leq \tau_j\}}$ converges weakly in
$L^2(\Omega\times [0,T]; \R^d))$ to the process
$\widehat{Z}_s^j1_{\{s\leq
\tau_j\}}$ as $k$ goes to infinity.\\
Now coming back to Equation (\ref{eqqq41}) we have for $k\geq n$
(and then $m_k^j \geq k\geq n$)

\begin{equation}\label{eqqq6}
\begin{array}{ll}
&2j\E\integ{0}{{\tau_j}}\psi'(Y_{s}^n -Y_{s}^{m_k^j})|Z_{s}^n-
Z_{s}^{m_k^j}|^2 ds +2j\E\integ{0}{{\tau_j}}|Z_{s}^n-
Z_{s}^{m_k^j}|^2 ds
\\&\leq \E\psi(Y_{{\tau_j}^-}^n-Y_{{\tau_j}^-}^{m_k^j})  +
\E\integ{0}{{\tau_j}}(\eta_s +j|Z_s^{m_k^j}|^2)\psi'(Y_{s}^n -Y_{s}^{m_k^j})ds\\
& \quad+\E\integ{0}{{\tau_j}} \psi'(Y_{s}^n -Y_{s}^{m_k^j})dA_s
+\E\displaystyle{\sum_{0<s< {\tau_j}}}\psi'(Y_{s-}^n
-Y_{s-}^{m_k^j})\, l_s
\\ & \leq \E\psi(Y_{{\tau_j}^-}^n-Y_{{\tau_j}^-}^{m_k^j})  +
\E\integ{0}{{\tau_j}}(\eta_s +2j|\widehat{Z}_s^{j}|^2)\psi'(Y_{s}^n
-Y_{s}^{m_k^j})ds \\
& \quad
+2j\E\integ{0}{{\tau_j}}|Z_s^{m_k^j}-\widehat{Z}_s^{j}|^2)\psi'(Y_{s}^n
-Y_{s}^{m_k^j})ds
\\ &\quad +\E\integ{0}{{\tau_j}} \psi'(Y_{s}^n
-Y_{s}^{m_k^j})dA_s +\E\displaystyle{\sum_{0<s<
{\tau_j}}}\psi'(Y_{s-}^n -Y_{s-}^{m_k^j})\, l_s.
\end{array}
\end{equation}
Since $|Z_{s}^n- Z_{s}^{m_k^j}|^2 = |Z_{s}^n- \widehat{Z}_{s}^{j}|^2
+2\langle \widehat{Z}_{s}^j- Z_{s}^{m_k^j},  Z_{s}^n-
\widehat{Z}_{s}^{j}\rangle + |\widehat{Z}_{s}^j- Z_{s}^{m_k^j}|^2$,
we have
\begin{equation}\label{eqqq6}
\begin{array}{ll}
&4j\E\integ{0}{{\tau_j}}\psi'(Y_{s}^n -Y_{s}^{m_k^j})\langle
\widehat{Z}_{s}^j- Z_{s}^{m_k^j},  Z_{s}^n-
\widehat{Z}_{s}^{j}\rangle ds +2j\E\integ{0}{{\tau_j}}|Z_{s}^n-
Z_{s}^{m_k^j}|^2 ds
\\ & \leq \E\psi(Y_{{\tau_j}^-}^n-Y_{{\tau_j}^-}^{m_k^j})  +
\E\integ{0}{{\tau_j}}(\eta_s +2j|\widehat{Z}_s^{j}|^2)\psi'(Y_{s}^n
-Y_{s}^{m_k^j})ds \\
&
 +\E\integ{0}{{\tau_j}} \psi'(Y_{s}^n
-Y_{s}^{m_k^j})dA_s +\E\displaystyle{\sum_{0<s<
{\tau_j}}}\psi'(Y_{s-}^n -Y_{s-}^{m_k^j})\, l_s.
\end{array}
\end{equation}
Letting $k$ to infinity, we get
 \begin{equation}\label{eqqq7}
\begin{array}{ll}
&2j\E\integ{0}{{\tau_j}}|Z_{s}^n- \widehat{Z}_{s}^{j}|^2 ds
\\ & \leq
2j\displaystyle{\liminf_{k\rightarrow
+\infty}}\E\integ{0}{{\tau_j}}|Z_{s}^n- Z_{s}^{m_k^j}|^2 ds
\\ & \leq \E\psi(Y_{{\tau_j}^-}^n-Y_{{\tau_j}}^{-})  + \E\integ{0}{{\tau_j}}(\eta_s
+2j|\widehat{Z}_s^{j}|^2)\psi'(Y_{s}^n
-Y_{s})ds \\
&
 \quad+\E\integ{0}{{\tau_j}} \psi'(Y_{s}^n
-Y_{s})dA_s +\E\displaystyle{\sum_{0<s< {\tau_j}}}\psi'(Y_{s-}^n
-Y_{s}^-)\, l_s.
\end{array}
\end{equation}
By dominated convergence theorem it follows that
$$\displaystyle\lim_{n\rightarrow +\infty}
\E\integ{0}{\tau_j}|Z_{s}^n- \widehat{Z}_{s}^{j}|^2 ds = 0.$$ 
By the uniqueness of the limit we obtain that $$
\widehat{Z}_s^j(\omega) 1_{\{0\leq s\leq \tau_j(\omega)\}} =
\widehat{Z}_s^{j+1}(\omega)1_{\{0\leq s\leq
\tau_j(\omega)\}},\,\,P(d\omega) ds-a.e.$$ For $s\in [0,T]$, let us
set $Z_s(\omega)= \displaystyle\lim_{j} \widehat{Z}_s^j 1_{\{ s \leq
\tau_j\}} = \widehat{Z}_s^{j(\omega)}(\omega),$ where $j(\omega)$ is
such that $\tau_{j(\omega)}(\omega) = T$. Then, for all $j\in\N$
$\E\integ{0}{\tau_j}\mid Z_s\mid^2 ds < +\infty$. Hence
$\integ{0}{T}\mid Z_s\mid^2 ds < +\infty, P-a.s. $ Moreover for all
$j\in\N$, we have
\begin{equation}\label{eqq10}
\displaystyle \lim_{n} \E\integ{0}{\tau_j}|Z_{s}^n -Z_{s}|^2 ds = 0.
\end{equation}
Proposition \ref{th1} is proved.\eop Let us now show that the
process $Y$ is \it{rcll}.
\begin{proposition}\label{pro00}
The process $Y$ is \it{rcll} and $Y_t^- = Y_{t-}$.
\end{proposition}
\bop. From Equation (\ref{eqqq41}) and according to
Bulkholder-Davis-Gundy inequality, there exists a universal constant
$C
>0$ such that
\begin{equation}\label{eqqq4}
\begin{array}{ll}
\E\displaystyle{\sup_{t<\tau_j}}\psi(Y_{t}^n -Y_{t}^{m}) &\leq
\E\psi(Y_{{\tau_j}^-}^n-Y_{{\tau_j}^-}^m)  +
\E\integ{0}{{\tau_j}}(\eta_s +j|Z_s^m|^2)\psi'(Y_{s}^n -Y_{s}^{m})ds\\
& +\E\integ{0}{{\tau_j}} \psi'(Y_{s}^n -Y_{s}^{m})dA_s
+\E\displaystyle{\sum_{0<s< {\tau_j}}}\psi'(Y_{s-}^n -Y_{s-}^{m})\,
l_s\\ & +C\E(\integ{0}{{\tau_j}}|\psi'(Y_{s}^n
-Y_{s}^{m})|^2|Z_{s}^n- Z_{s}^{m}|^2 ds)^{\frac12}.
\end{array}
\end{equation}
Hence $$ \lim_{n} \E\sup_{s< \tau_j}\psi(Y_{s}^n -Y_{s})=0.$$ It
follows that $$ \lim_{n} \E\sup_{s< \tau_j}\mid Y_{s}^n -Y_{s}\mid
=0.$$  Consequently  $Y$ is \it{rcll} and $Y_t^- = Y_{t-}$, since
$P[\displaystyle{\cup_{j=1}^{n}}(\tau_j= T)] =1$.\eop
\begin{proposition}
For all $j\in\N^*$, we have
\begin{enumerate}
\item $  \displaystyle{\lim_{n\rightarrow +\infty}}
\E\integ{0}{\tau_j}\mid f_{n}(s, Y_s^{n}, Z_s^{n}) - f(s, Y_s,
Z_s)\mid ds =0.$
\item $\displaystyle{\lim_{n\rightarrow +\infty}}\E\integ{0}{\tau_j}\mid g_{n}(s,
Y_s^{n}) - g(s, Y_s)\mid dA_s =0.$
\item $
\displaystyle{\lim_{n\rightarrow +\infty}}\E\displaystyle{\sum_{0<s<
{\tau_j}}}|h_n(s, Y_{s-}^{n},Y_s^n)-h(s,Y_{s-}, Y_s)| =0. $
\end{enumerate}
\end{proposition}

\bop. In view of (\ref{eqq10}) there exists a subsequence
$\bigg(n_k^j\bigg)_k$ of $n$ such that :
\begin{enumerate}
\item
$ \E\integ{0}{\tau_j}|Z_{s}^{n_k^j} -Z_{s}|^2 ds\leq
\frac{1}{2^k}\quad \mbox{and then} \quad
\E\integ{0}{\tau_j}\sum_{k=0}^{+\infty}|Z_{s}^{n_k^j} -Z_{s}|^2
ds\leq 2, $
\item $Z_s^{n_k^j}(\omega)\longrightarrow Z_s(\omega),\, a.e.\,\,
(s,\omega)\in[0,\tau_j]\times \Omega,\quad\mbox{and}\quad \mid
Z_s^{n_k^j}(\omega)\mid \leq h^j_s,\,\, a.e.\,\,
(s,\omega)\in[0,\tau_j]\times \Omega,$ where $h_s^j = 1_{\{ s \leq
\tau_j\}}
\bigg(2|Z_s|^2+2\displaystyle\sum_{k=0}^{+\infty}|Z_s^{n_k^j}-Z_s|^2\bigg)^{\frac12}$.
\end{enumerate}
 It follows then from Lemma \ref{lem0} that
$$
\begin{array}{ll}
& \E\integ{0}{\tau_j}| f_{n_k^j}(s, Y_s^{n_k^j}, Z_s^{n_k^j}) - f(s,
Y_s^{n_k^j}, Z_s^{n_k^j})| ds \\ &= \E\integ{0}{\tau_j} f_{n_k^j}(s,
Y_s^{n_k^j}, Z_s^{n_k^j}) - f(s, Y_s^{n_k^j}, Z_s^{n_k^j}) ds  \\ &
= \E\integ{0}{\tau_j}\bigg(f_{n_k^j}(s, Y_s^{n_k^j}, Z_s^{n_k^j}) -
f(s, Y_s^{n_k^j}, Z_s^{n_k^j})\bigg)
1_{\{ \mid Z_s^{n_k^j} -Z_s\mid \leq 1\}}ds\\
&\quad +\E\integ{0}{\tau_j}\bigg(f_{n_k^j}(s, Y_s^{n_k^j},
Z_s^{n_k^j}) - f(s, Y_s^{n_k^j}, Z_s^{n_k^j})\bigg) 1_{\{ \mid
Z_s^{n_k^j} -Z_s\mid \geq 1\}}ds
\\ &\leq  \E\integ{0}{\tau_j}\sup_{(y,z)\in [0,1]\times B(Z_s,1)}\bigg(f_{n_k^j}(s, y, z) -
f(s, y,z)\bigg) ds -\E\integ{0}{\tau_j} f(s, Y_s^{n_k^j},
Z_s^{n_k^j}) 1_{\{ \mid Z_s^{n_k^j} -Z_s\mid \geq 1\}}ds
\\ &\leq
\E\integ{0}{\tau_j}\sup_{(y,z)\in [0,1]\times
B(Z_s,1)}\bigg(f_{n_k^j}(s, y, z) - f(s, y,z)\bigg) ds
\\ &
\quad +\E\integ{0}{\tau_j}\bigg({\eta}_s +2j \mid Z_s^{n_k^j}
-Z_s\mid^2 + 2j \mid Z_s\mid^2 \bigg)(\mid Z_s^{n_k^j}
-Z_s\mid\wedge 1)ds,
\end{array}
$$
where $B(Z,1)$ is the closed ball of center $Z$ and radius $1$.\\
By taking account of Lemma \ref{lem0} one can see that
$$\displaystyle{\sup_{(y,z)\in [0,1]\times
B(Z_s,1)}\bigg(f_{n_k^j}(s, y, z) - f(s, y,z)\bigg)},$$ converges
pointwise. But, on $[0, \tau_j[$, we have also that
$$
\sup_{(y,z)\in [0,1]\times B(Z_s,1)}\bigg(f_{n_k^j}(s, y, z) - f(s,
y,z)\bigg) \leq \eta_s +j(|Z_s| +1)^2.
$$
Henceforth, by using Lebesgue's dominated convergence theorem, we
get
$$ \lim_k\E\integ{0}{\tau_j}\sup_{(y,z)\in [0,1]\times
B(Z_s,1)}\bigg(f_{n_k^j}(s, y, z) - f(s, y,z)\bigg) ds = 0, $$ and
$$\lim_k\E\integ{0}{\tau_j}\bigg({\eta}_s +j \mid Z_s^{n_k^j} -Z_s\mid^2 +
j \mid Z_s\mid^2 \bigg)(\mid Z_s^{n_k^j} -Z_s\mid\wedge 1)ds = 0.
$$
Therefore
$$ \lim_k \E\integ{0}{\tau_j}\bigg(f_{n_k^j}(s, Y_s^{n_k^j},
Z_s^{n_k^j}) - f(s, Y_s^{n_k^j}, Z_s^{n_k^j})\bigg) ds =0.$$ It
follows also from Lebesgue's dominated convergence theorem and the
continuity of $f$ that for all $j\in\N$
 $$  \lim_k \E\integ{0}{\tau_j}\mid f(s, Y_s^{n_k^j},
Z_s^{n_k^j}) - f(s, Y_s, Z_s)\mid ds =0.$$ Hence for all $j\in\N$ $$
\lim_k \E\integ{0}{\tau_j}\mid f_{n_k^j}(s, Y_s^{n_k^j},
Z_s^{n_k^j}) - f(s, Y_s, Z_s)\mid ds =0.$$ Since the above limit
doesn't depend on the choice of the subsequence $(n_k^j)_{k}$ we
have for all $j\in\N$
$$  \lim_n
\E\integ{0}{\tau_j}\mid f_{n}(s, Y_s^{n}, Z_s^{n}) - f(s, Y_s,
Z_s)\mid ds =0.$$ It not difficult also to prove that for all
$j\in\N$
$$  \lim_n
\E\integ{0}{\tau_j}\mid g_{n}(s, Y_s^{n}) - g(s, Y_s)\mid dA_s =0.
$$
Now
$$
\begin{array}{ll}
&\E\displaystyle{\sum_{0<s< {\tau_j}}}|h_n(s,
Y_{s-}^{n},Y_s^n)-h(s,Y_{s-}, Y_s)|
\\ & \leq
\E\displaystyle{\sum_{i=1}^{n}}|h({T_i},
Y_{{T_i}^-}^{n},Y_{T_i}^n)-h(s,Y_{{T_i}^-},
Y_{T_i})|1_{\{T_i<\tau_j\}} +
\E\displaystyle{\sum_{i=n+1}^{+\infty}}l_{T_i}1_{\{T_i<\tau_j\}}
\end{array}
$$
Since $|h({T_i}, Y_{{T_i}^-}^{n},Y_{T_i}^n)-h(s,Y_{{T_i}^-},
Y_{T_i})|\leq 2l_{T_i}$ and
$\E\displaystyle{\sum_{i=1}^{+\infty}}l_{T_i}1_{\{T_i<\tau_j\}}\leq
j$, by Lebesgue's convergence theorem and the continuity of $h$ it
follows that
$$
\lim_{n\rightarrow +\infty}\E\displaystyle{\sum_{0<s<
{\tau_j}}}|h_n(s, Y_{s-}^{n},Y_s^n)-h(s,Y_{s-}, Y_s)| =0.
$$
 \eop
\subsubsection{Identification of the limit}
\begin{proposition}\label{pro000} The process $(Y, Z, K^+, K^-)$ defined in Subsection \ref{sub4.4.2} satisfies,
$P$-a.s., the following:
\begin{enumerate}
\item $K_{T}^{+}< +\infty$.
\item
$$
\begin{array}{ll}
{Y}_{t} & = \integ{t}{T}{f}(s,{Y}_{s},{Z}_{s})ds +
 \integ{t}{T}g(s,{Y}_{s})d{A}_s + \integ{t}{T}
 d{R}_s
 +\displaystyle{\sum_{t<s\leq
T}}h(s,Y_{s-}, Y_s)+\integ{t}{T}d{K}_{s}^{+}
\\ &\qquad -\integ{t}{T}d{K}_s^{-}
-\integ{t}{T}{Z}_{s}dB_{s}.
\end{array}
$$
\item $\integ{0}{T}( {U}_{t-}-Y_{t-}) dK_{t}^{-} = \integ{0}{T}(Y_{t-} -L_{t-}) dK_{t}^{+} =
0$ and $dK^+ \perp dK^-$.
\end{enumerate}
\end{proposition}
\bop. \it{1.} From Equation (\ref{eq9})$(i)$ we obtain, $\forall
j$,\, $
 \sup_{n}\E K_{\tau_j}^{n+} <+\infty$. It follows then from Fatou's lemma that for any $j\in\N$,\, $ \E
K_{\tau_{j-}}^{+}<+\infty$. Henceforth $K_{T-}^{+}< +\infty, P$-a.s.
But, $P$-a.s.
$$
\begin{array}{ll}
\Delta_T K^{+} & = \displaystyle{\lim_{n\rightarrow
+\infty}}\Delta_T K^{n+}
\\ &= \displaystyle{\lim_{n\rightarrow +\infty}}(L_{T-} - (h_n(T,
Y_{T-}^n, 0)+\Delta_T R))^+
\\ & = (L_{T-} - (h(T,
Y_{T-}, 0)+\Delta_T R))^+ \\ & < +\infty,
\end{array}
$$
then $K_{T}^{+}< +\infty, P$-a.s.\\
\it{2.} Since
$$
\begin{array}{ll}
{Y}_{t}^n & ={Y}_{{\tau_{j}}^-}^n
+\integ{t}{\tau_{j}}{f}_n(s,{Y}_{s}^n,{Z}_{s}^n)ds +
 \integ{t}{\tau_{j}}g_n(s,{Y}_{s}^n)d{A}_s + \integ{t}{{\tau_{j}}^-}
 d{R}_s
 \\ &+\displaystyle{\sum_{t<s<
\tau_j}}h_n(s,Y_{s-}^n, Y_s^n)+\integ{t}{{\tau_{j}}^-}d{K}_{s}^{n+}
 -\integ{t}{{\tau_{j}}^-}d{K}_s^{n-}
-\integ{t}{\tau_j}{Z}_{s}^ndB_{s}
\end{array}
$$
Passing to the limit as $n$ goes to infinity and using the fact that
$\tau_j$ is a stationary stopping time we get $P-$a.s
$$
\begin{array}{ll}
{Y}_{t} & ={Y}_{T-} +\integ{t}{T}{f}(s,{Y}_{s},{Z}_{s})ds +
 \integ{t}{T}g(s,{Y}_{s})d{A}_s + \integ{t}{T-}
 d{R}_s
 \\ &+\displaystyle{\sum_{t<s<
T}}h(s,Y_{s-}, Y_s)+\integ{t}{T-}d{K}_{s}^{+}
 -\integ{t}{T-}d{K}_s^{-}
-\integ{t}{T}{Z}_{s}dB_{s}
\\ & = \integ{t}{T}{f}(s,{Y}_{s},{Z}_{s})ds +
 \integ{t}{T}g(s,{Y}_{s})d{A}_s + \integ{t}{T}
 d{R}_s
 +\displaystyle{\sum_{t<s\leq
T}}h(s,Y_{s-}, Y_s)+\integ{t}{T}d{K}_{s}^{+}
\\ & -\integ{t}{T}d{K}_s^{-}
-\integ{t}{T}{Z}_{s}dB_{s} +(Y_{T-}  -h(T, Y_{T-}, 0)-\Delta_T
K^++\Delta_T K^-).
\end{array}
$$
By using relation (\ref{equa1}) we obtain
$$(Y_{T-} -\xi -h(T, Y_{T-}, 0)-\Delta_T K^++\Delta_T
K^-)=\displaystyle{\lim_{n\rightarrow +\infty}}(Y_{T-}^n -\xi-h_n(T,
Y_{T-}^n, 0) -\Delta_T K^{n+}+\Delta_T K^{n-})=0,
$$
then it follows that, $P-$a.s.
$$
\begin{array}{ll}
{Y}_{t} & = \integ{t}{T}{f}(s,{Y}_{s},{Z}_{s})ds +
 \integ{t}{T}g(s,{Y}_{s})d{A}_s + \integ{t}{T}
 d{R}_s
 +\displaystyle{\sum_{t<s\leq
T}}h(s,Y_{s-}, Y_s)+\integ{t}{T}d{K}_{s}^{+}
 \\ &\qquad -\integ{t}{T}d{K}_s^{-}
-\integ{t}{T}{Z}_{s}dB_{s}.
\end{array}
$$
\it{3.} Now let us prove the minimality conditions. We have
$$
\integ{0}{T}( {U}_{t-}-Y_{t-}^{n}) dK_{t}^{n-}=0.
$$
Hence, since $dK^{-} = \displaystyle\inf_{n}dK^{n-}$, we get
$$
\integ{0}{T}( {U}_{t-}-Y_{t-}^{n}) dK_{t}^{-}=0.
$$
 It follows then from Fatou's lemma that
$$
\integ{0}{T}( {U}_{t-}-Y_{t-}) dK_{t}^{-}  = 0.
$$
On the other hand
$$
\integ{0}{T}(Y_{t-}^{n} -L_{t-}) dK_{t}^{n+}  = 0.
$$
Hence, since $Y = \displaystyle\inf_{n}Y^{n}$, we obtain
$$
\integ{0}{T}(Y_{t-} -L_{t-}) dK_{t}^{n+} = 0.
$$
Applying Fatou's lemma we obtain
$$
\integ{0}{T}(Y_{t-} -L_{t-}) dK_{t}^{+} = 0.
$$
Now, since $dK^{+} = \displaystyle\sup_{n}dK^{n+}$,\,\, $dK^{-} =
\displaystyle\inf_{n}dK^{n-}$ and the measures $dK^{n+}$ and
$dK^{n-}$ are singular, it follows that $dK^{+}$ and $dK^{-}$ are
singular. \eop

\bop\,\,\bf{of Theorem \ref{the122}}.  By using Propositions
\ref{pro00}-\ref{pro000}, it is not difficult to see that the
process $(Y, Z, K^+ , K^-)$ satisfies Equation (\ref{eq221}). It
remains to prove $(Y, Z, K^+ , K^-)$ is maximal.\\  Let $(Y', Z',
K'^{+} , K'^{-})$ be an another solution to Equation (\ref{eq221}).
By comparison theorem (Theorem \ref{th111}) we have that $Y' \leq
Y^n$ and then $Y'\leq Y$. The proof of Theorem \ref{the122} is then
finished.\eop
\section{Comparison theorem for maximal solutions }
 This section is devoted to show a comparison theorem for
maximal solutions. For this reason, suppose that assumptions
$(\bf{A.1})$--$(\bf{A.4})$ hold and $(Y, Z, K^+, K^-)$ is the
maximal solution for the following GRBSDE
\begin{equation}
\label{eq0q} \left\{
\begin{array}{ll}
(i) & 
 Y_{t}=\xi
+\integ{t}{T}f(s,Y_{s},Z_{s})ds+\dint_t^Tg(s,
Y_s)dA_s +\sum_{t<s\leq T}h(s, Y_{s-}, Y_s)\\
&\qquad\quad+\integ{t}{T}dK_{s}^+ -\integ{t}{T}dK_{s}^-
-\integ{t}{T}Z_{s}dB_{s}\,, t\leq T,
\\ (ii)&
\forall t\in[0,T[,\,\, L_t \leq Y_{t}\leq U_{t},\\  (iii)&
\integ{0}{T}( Y_{t-}-L_{t-})
dK_{t}^+= \integ{0}{T}( U_{t-}-Y_{t-}) dK_{t}^-=0,\,\, \mbox{a.s.}, \\
(iv)& Y\in {\cal D}, \quad K^+, K^-\in {\cal K}, \quad Z\in {\cal
L}^{2,d},  \\ (v)& dK^+\perp  dK^-.
\end{array}
\right. \end{equation} Let $(Y', Z', K'^+, K'^-)$ be a solution for
the following GRBSDE
\begin{equation}
\label{eq1q} \left\{
\begin{array}{ll}
(i) & 
 Y'_{t}=\xi'
+\integ{t}{T}f'(s)ds+\dint_t^Tg'(s)dA_s +\sum_{t<s\leq
T}h'(s)+\integ{t}{T}dK'^{+}_{s}\\
&\qquad\quad -\integ{t}{T}dK'^{-}_{s} -\integ{t}{T}Z'_{s}dB_{s}\,,
t\leq T,
\\ (ii)&
\forall t\in[0,T[,\,\, L'_t \leq Y'_{t}\leq U'_{t},\\  (iii)&
\integ{0}{T}( Y'_{t-}-L'_{t-})
dK'^{+}_{t}= \integ{0}{T}( U'_{t-}-Y'_{t-}) dK'^{-}_{t}=0,\,\, \mbox{a.s.}, \\
(iv)& Y'\in {\cal D}, \quad K'^+, K'^-\in {\cal K}, \quad Z'\in
{\cal L}^{2,d},  \\ (v)& dK'^+\perp  dK'^-.

\end{array}
\right. \end{equation} Assume moreover that for every $t\in[0, T]$
\begin{enumerate}
\item $\xi\leq \xi'$,
\item $Y'_t\leq U_t$,
\item $L'_t\leq Y_t$,
\item $f'(s)ds\leq f(s, Y'_s, Z'_s)ds$ on $[0, T]$,
\item $g'(s)dA_s\leq g(s, Y'_s)dA_s$ on $[0, T]$,
\item $h'(s)\leq h(s, Y'_{s-}, Y'_s)$ for every $s\in]0, T]$.
\end{enumerate}
The way in which the maximal solution for GRBSDE (\ref{eq0}) has
been constructed allow us to deduce the following comparison
theorem.
\begin{theorem} (Comparison theorem for maximal solutions)
\label{th103} Assume that the above assumptions hold then we have :
\begin{enumerate}
\item $Y'_t\leq Y_t$, for $t\in [0,T]$, $P-$a.s.
\item $
1_{\{U'_{t-} = U_{t-}\}}dK'^{-}_t \leq dK^{-}_t\,\,\,
\mbox{and}\,\,\, 1_{\{L'_{t-} = L_{t-}\}}dK^{+}_t \leq dK'^{+}_t. $
\end{enumerate}
\end{theorem}
\bop. By using an exponential change like the one used in Section
\ref{sec3.2}, one can suppose assumptions $(\bf{H.0})$--$(\bf{H.4})$
instead of $(\bf{A.1})$--$(\bf{A.4})$. Moreover by using the
approximations scheme of Section \ref{sec4.4} and defining $(Y^n,
Z^n, K^{n+}, K^{n-})$ as the unique solution of GRBSDE (\ref{eq9})
it follows from comparison theorem (Theorem \ref{th111}) that for
every $n\in\N$
\begin{enumerate}
\item $Y'_t\leq Y_t^n$, for $t\in [0,T]$, $P-$a.s.
\item $
1_{\{U'_{t-} = U_{t-}\}}dK'^{-}_t \leq dK^{n-}_t\,\,\,
\mbox{and}\,\,\, 1_{\{L'_{t-} = L_{t-}\}}dK^{n+}_t \leq dK'^{+}_t. $
\end{enumerate}
By letting $n$ to infinity, the result follows by using the
convergence result of Section \ref{sec4.4}. Theorem \ref{th103} is
proved. \eop
\begin{remark} It should be noted that the result of Theorem \ref{th103} remains true
if the data $f', g'$ and $h'$ of GRBSDE (\ref{eq1q}) depend on $(Y',
Z')$, $Y'$ and $(Y'_{-}, Y')$ respectively.
\end{remark}
\section{\Large{\bf{Appendix : Comparison theorem}}}
The comparison theorem for real-valued BSDEs turns out to be one of
the classic results of the theory of BSDE. It allows to compare the
solutions of two real-valued BSDEs whenever we can compare the
terminal conditions and the generators. This section is devoted to
show a comparison theorem for the following GRBSDE :
\begin{equation}
\label{eq000000} \left\{
\begin{array}{ll}
(i) & 
 Y_{t}=\xi
+\integ{t}{T}f(s,Y_{s},Z_{s})ds+\dint_t^Tg(s,
Y_s)dA_s +\dint_t^T dR_s+\sum_{t<s\leq T}h(s, Y_{s-}, Y_s)\\
&\qquad\quad+\integ{t}{T}dK_{s}^+ -\integ{t}{T}dK_{s}^-
-\integ{t}{T}Z_{s}dB_{s}\,, t< T,
\\ (ii)&
\forall t\in[0,T[,\,\, L_t \leq Y_{t}\leq U_{t},\\  (iii)&
\integ{0}{T}( Y_{t-}-L_{t-})
dK_{t}^+= \integ{0}{T}( U_{t-}-Y_{t-}) dK_{t}^-=0,\,\, \mbox{a.s.}, \\
(iv)& Y\in {\cal D}, \quad K^+, K^-\in {\cal K}, \quad Z\in {\cal
L}^{2,d},  \\ (v)& dK^+\perp  dK^-.
\end{array}
\right. \end{equation} Let $(Y^i, Z^i, K^{i+}, K^{i+})$ $(i=1,2)$ be
two solutions of Equation (\ref{eq000000}) associated respectively
with \\ $(\xi^1, f^1, g^1, h^1, A^1, L^1, U^1)$ and $(\xi^2, f^2,
g^2, h^2, A^2, L^2, U^2)$, such
that, for $(i=1,2)$ the following assumptions are satisfied :\\
\noindent$(\bf{D.1})$  $L^i: [0,T]\times \Omega\longrightarrow \R$
and  $U^i: [0,T]\times \Omega\longrightarrow \R$ are two \it{rcll}
barriers processes satisfying
$$
L^1_t\leq Y^2_t, \quad Y^1_t\leq U^2_t, \quad U^1_{t}\wedge
U^2_t-L^1_t\vee L^2_{t} \leq 2\quad\forall t\in [0,T[,\,\, P-a.s.
$$
\noindent$(\bf{D.2})$ $0\leq A_T^2\leq 1$, $\quad dR_t^1\leq dR_t^2$ on $[0, T]$.\\
 \noindent$(\bf{D.3})$ $f^1$ and $f^2$ are such that : \begin{itemize} \item[\bf a.] $ f^1(s,{Y}_{s}^1,{Z}_{s}^1)\leq
f^2(s,{Y}_{s}^1,{Z}_{s}^1)$,\,\,\, $dsP(d\omega)-a.e.$.
\item[\bf b.] $f^2$ is uniformly Lipschitz with respect to $(y,z)$
with Lipschitz constant $C_1\geq 0$.
\end{itemize}
\noindent$(\bf{D.4})$ $g^1$ and $g^2$ are such that :
\begin{itemize}
\item[\bf a.] $ g^1(s,{Y}_{s}^1)dA_s^1\leq g^2(s,{Y}_{s}^1)dA_s^2$
on $[0, T],\,$ $P(d\omega)-$a.s.
\item[\bf b.] $g^2$ is uniformly Lipschitz with respect to $y$
with Lipschitz constant $C_2\geq 0$.
\end{itemize}
\noindent$(\bf{D.5})$ $h^1$ and $h^2$ are such that :
\begin{itemize}
\item[\bf a.]$P$-a.s.,\,\,$\forall s\leq T$,\,\,
$h^1(s,Y^1_{s-}, Y^1_s)\leq h^2(s, Y^1_{s-}, Y^1_s)$,
\item[\bf b.] $P$-a.s., $\forall (t, x)\in]0, T]\times \R$,\,\, the function $y\mapsto y+h^2(t,
\omega, L_{t-}^2(\omega)\vee x \wedge U_{t-}^2(\omega),
L_t^2(\omega)\vee y \wedge U_t^2(\omega))$ is nondecreasing.
\item[\bf c.] there exists a family of stopping times $S_0 =0\leq S_1\leq ...\leq S_{p+1}=T$ such that : $\forall s\notin \{S_1,...,S_p, S_{p+1}=T\}$
$h^2(s,x, y)=0,\,\,$ for every $(x,y)\in\R^2$ and for each
$i\in\{1,...,p\}$
$$
Y_{S_{i-}}^2 = \max\{ x\in[L_{S_{i-}}^2, U_{S_{i-}}^2]\,\, : \,\, x=
L_{S_{i-}}^2\vee[ Y_{S_{i}}^2+ h^2(t,x,Y_{S_{i}}^2)+\Delta_{S_{i}}
R^2]\wedge U_{S_{i-}}^2 \}.
$$
\end{itemize}

The following comparison theorem plays a crucial role in our proofs.
\begin{theorem}\label{th111} (Comparison theorem)  Assume that assumptions $(\bf{D.1})-(\bf{D.5})$
hold. Then
\begin{enumerate}
\item $Y^{1}_t\leq Y^{2}_ t$, for $t\in [0,T]$, $P-$a.s.
\item $
1_{\{U^1_{t-} = U^2_{t-}\}}dK^{1-}_t \leq dK^{2-}_t\,\,\,
\mbox{and}\,\,\, 1_{\{L^1_{t-} = L^2_{t-}\}}dK^{2+}_t \leq
dK^{1+}_t. $
\end{enumerate}
\end{theorem}
In order to prove Theorem \ref{th111}. we need the following lemma.
\begin{lemma}\label{lem4}
let $\tau\in [0,T]$ be a stopping time. If $Y^{1}_{\tau}\leq
Y^{2}_{\tau}$ then $Y^{1}_{\tau -}\leq Y^{2}_{\tau -}.$
\end{lemma}
\bop. We distinguish two cases. \\
\ni 1. If $Y_{{\tau -}}^1\leq L^1_{\tau -} \vee L^2_{\tau -}$ then
it is obvious, by assumption $(\bf{D.1})$, that $Y^{1}_{\tau -}\leq
Y^{2}_{\tau -}$.
\\
\ni 2. If $Y_{{\tau -}}^1 > L^1_{\tau -} \vee L^2_{\tau -}$, by
Lemma \ref{lem1} and assumptions $(\bf{D.2})-(\bf{D.5})$ , we have
$$
\begin{array}{lll}
Y_{{\tau -}}^1 &=  [Y_{\tau}^1 + h^1(\tau, Y_{{\tau -}}^1,
Y_{\tau}^1)+\Delta_{t} R^1]\wedge U_{{\tau -}}^1
\\ & \leq  [Y_{\tau}^2 + h^2(\tau,Y_{{\tau -}}^1, Y_{\tau}^2)+\Delta_{t} R^2]\wedge U_{{\tau -}}^2.
\end{array}
$$
Since $Y_{t-}^2 = \max\{ x\in[L_{t-}^2, U_{t-}^2]\,\, : \,\, x=
L_{t-}^2\vee[ Y_{t-}^2+ h^2(t,x,Y_{t-}^2)+\Delta_{t} R^2]\wedge
U_{t-}^2 \}$, then $Y^{1}_{\tau -}\leq Y^{2}_{\tau -}$.\eop

\bop\,\, \bf{of Theorem \ref{th111}.}\,\,\it{1.} We proceed by
induction. We have $Y_{S_{p+1}}^1 =\xi^1\leq Y_{S_{p+1}}^2=\xi^2$.
Suppose that for $i\in\{0,...p \}$,\,\,
$Y_{S_{i+1}}^1\leq Y_{S_{i+1}}^2\,\,P-$a.s., 
then by Lemma \ref{lem4}. we have $Y_{S_{i+1}^-}^1\leq
Y_{S_{i+1}^-}^2$. Let $\tau\in [S_{i}, S_{i+1}]$ and define
$$
\lambda_{\tau}:=\inf\{s > \tau :Y^2_{s-}\geq Y^1_{s-}\}\wedge
S_{i+1}.
$$ 
On the set $\{w\in\Omega : \tau(\omega)< \lambda_{\tau}(\omega)\}$,
for every $s\in ]\tau, \lambda_{\tau}[$ we have $Y^2_{s-}< Y^1_{s-}$
and then for every $s\in [\tau, \lambda_{\tau}[$,\,\, $Y^2_s\leq
Y^1_s$, hence $Y^2_{\lambda_{\tau}^-}\leq Y^1_{\lambda_{\tau}^-}$.
Now if $Y^2_{\lambda_{\tau}^-}< Y^1_{\lambda_{\tau}^-}$, then there
exists a sequence $(v^n)_{n\geq 1} \in ]\lambda_{\tau}, S_{i+1}[$
which converges to $\lambda_{\tau}$ such that $Y^2_{v_{n-}}\geq
Y^1_{v_{n-}}$. By letting $n$ to infinity we have
$Y^2_{\lambda_{\tau}}\geq Y^1_{\lambda_{\tau}}$. By Lemma
\ref{lem4}. it follows that $Y^2_{\lambda_{\tau}^-}\geq
Y^1_{\lambda_{\tau}^-}$, which is absurd. Then
$Y^2_{\lambda_{\tau}^-} = Y^1_{\lambda_{\tau}^-}$.


For all $t\in [\tau, \lambda_{\tau}[$ we have
$$
\begin{array}{ll}
&Y_{t}^1 -Y_{t}^{2} \\&=
\integ{t}{\lambda_{\tau}}(f^1(s,Y_{s}^1,Z_{s}^1)
-f^2(s,Y_{s}^{2},Z_{s}^{2}))ds +
\integ{t}{\lambda_{\tau}}g^1(s,Y_{s}^1)dA_s^1
-\integ{t}{\lambda_{\tau}}g^2(s,Y_{s}^{2}))dA_s^2
\\ &+\integ{t}{\lambda_{\tau -}}(dR^1_s-dR^2_s) +\displaystyle{\sum_{t< s<
\lambda_{\tau}}}(h^1(s, Y_{s-}^1, Y_s^1)-h^2(s, Y_{s-}^2,
Y_s^2))
\\
&+\integ{t}{\lambda_{\tau -}} d(K_{s}^{1+}-K_{s}^{2+})
-\integ{t}{\lambda_{\tau -}}d(K_{s}^{1-}-K_{s}^{2-})
-\integ{t}{\lambda_{\tau}}(Z_{s}^1 - Z_{s}^{2})dB_{s}.
\end{array}
$$
But for every $t\in [\tau, \lambda_{\tau}[$,\,\,
$\integ{t}{\lambda_{\tau -}}dK_{s}^{1+} = \integ{t}{\lambda_{\tau
-}}1_{\{s < \lambda_{\tau}\}} 1_{\{Y_{s-}^2 <Y_{s-}^1 =
L_{s-}^1\}}dK_{s}^{1+} = 0$ and by the same way
$\integ{t}{\lambda_{\tau -}}dK_{s}^{2-} = 0$. 
Therefore
$$
\begin{array}{ll}
&Y_{t}^1 -Y_{t}^{2} \\&=
\integ{t}{\lambda_{\tau}}(f^1(s,Y_{s}^1,Z_{s}^1)
-f^2(s,Y_{s}^{2},Z_{s}^{2}))ds +
\integ{t}{\lambda_{\tau}}g^1(s,Y_{s}^1)dA_s^1
-\integ{t}{\lambda_{\tau}}g^2(s,Y_{s}^{2}))dA_s^2\\
&+\integ{t}{\lambda_{\tau -}}(dR^1_s-dR^2_s) +\displaystyle{\sum_{t<
s< \lambda_{\tau}}}(h^1(s, Y_{s-}^1, Y_s^1)-h^2(s, Y_{s-}^2,
Y_s^2))\\ & -\integ{t}{\lambda_{\tau -}} d(K_{s}^{1-}+K_{s}^{2+})
-\integ{t}{\lambda_{\tau}}(Z_{s}^1 - Z_{s}^{2})dB_{s}.
\end{array}
$$
Since $f^2$ and $g^2$ are Lipschitz then we can write
$f^2(s,Y_{s}^{1},Z_{s}^{1})- f^2(s,Y_{s}^{2},Z_{s}^{2}) = a_s (Y_s^1
-Y_s^2 ) + \langle \widehat{b}_s, Z_s^1 - Z_s^2 \rangle,\,\,s\leq T$
and $g^2(s,Y_{s}^{1})- g^2(s,Y_{s}^{2}) = \widehat{a}_s (Y_s^1
-Y_s^2 )$, where $(a_t)_{t\leq T}$, $(\widehat{a}_t)_{t\leq T}$ and
$(\widehat{b}_t)_{t\leq T}$ are bounded ${\cal P}$-measurable
processes, we have
$$
\begin{array}{ll}
&Y_{t}^1 -Y_{t}^{2} \\&= \integ{t}{\lambda_{\tau}}\bigg(a_s(Y^1_s
-Y^2_s) + \langle \widehat{b}_s, Z^1_s -Z^2_s\rangle\bigg)ds +
\integ{t}{\lambda_{\tau}}\widehat{a}_s(Y^1_s -Y^2_s)dA_s^2
\\ &-\integ{t}{\lambda_{\tau}}1_{\{s < \lambda_{\tau}\}} dV_s
-\integ{t}{\lambda_{\tau}}(Z_{s}^1 - Z_{s}^{2})dB_{s},
\end{array}
$$
where the process $V\in {\cal K}$ is defined by :
$$
\begin{array}{ll}
V_t := &\integ{0}{t} (f^2(s,Y_{s}^{1},Z_{s}^{1})-
f^1(s,Y_{s}^{1},Z_{s}^{1}))ds + \integ{0}{t} g^2(s,Y_{s}^{1})
dA^2_s- \integ{0}{t}g^1(s,Y_{s}^{1}))dA^1_s\\ &  +
\displaystyle{\sum_{0< s< t}}(h^2(s, Y_{s-}^1, Y_{s}^1)-h^1(s,
Y_{s-}^1, Y_s^1))+ R^2_t -R^1_t+ K_t^{2+}+ K_t^{1-}.
\end{array}
$$
Setting $\Gamma_t = e^{ \int_{\tau}^t (a_s ds + \widehat{a}_s
dA_s^2)+\int_{\tau}^t\widehat{b}_s dB_s -\frac12 \int_{\tau}^t\mid
\widehat{b}_s\mid^2 ds}$, it follows that
$$
\begin{array}{ll}
\Gamma_t(Y_{t}^1-Y_{t}^2) = -\integ{t}{\lambda_{\tau}}1_{\{s <
\lambda_{\tau}\}} \Gamma_s dV_s -\integ{t}{\lambda_{\tau}}\Gamma_s
(Z_{s}^1 - Z_{s}^{2})dB_{s}.
\end{array}
$$
Let $\theta_n$ be a family of stopping times defined by
$$
\theta_n =\inf\{s\geq \tau: \int_{\tau}^{s} \Gamma_{s}|Z_{s}^1 -
Z_{s}|^2ds \geq n\}\wedge \lambda_{\tau}.
$$
By assumption $(\bf{D.1})$ we have
$$
\begin{array}{ll}
\E 1_{\{\tau < \lambda_{\tau}\}}(Y_{\tau}^1-Y_{\tau}^2) &\leq \E
1_{\{\theta_n <
\lambda_{\tau}\}}\Gamma_{\theta_n}(Y_{\theta_{n}}^1-Y_{\theta_{n}}^2)
\\ &\leq \E 1_{\{\theta_n <
\lambda_{\tau}\}}\Gamma_{\lambda_{\tau}}(U_{\theta_{n}}^1 \wedge
U_{\theta_{n}}^2-L_{\theta_{n}}^1\vee L_{\theta_{n}}^2)^+
\\ & \leq 2\E 1_{\{\theta_n <
\lambda_{\tau}\}}\Gamma_{\lambda_{\tau}}.
\end{array}
$$
Since $P[\displaystyle{\cup_{n\geq 0}}(\theta_n =\lambda_{\tau})]
=1$ and $\E \Gamma_{\lambda_{\tau}} <+\infty$, it follows that
$$
\begin{array}{ll}
\E 1_{\{\tau < \lambda_{\tau}\}}(Y_{\tau}^1-Y_{\tau}^2)
 &\leq 2\displaystyle{\lim_{n\rightarrow +\infty}}\E 1_{\{\theta_n <
\lambda_{\tau}\}}\Gamma_{\lambda_{\tau}}
\\ & =2\E 1_{\{\cap_{n\geq 1}\{\theta_n <
\lambda_{\tau}\}\}}\Gamma_{\lambda_{\tau}}\\ & =0.
\end{array}
$$
 Hence $1_{\{\tau < \lambda_{\tau}\}}(Y_{\tau}^1-Y_{\tau}^2) = 0,
 \,\,P-$a.s. and then $P[\tau<\lambda_{\tau},\,Y_{\tau}^2<Y_{\tau}^1 ] =0$. Hence
 $\tau = \lambda_{\tau}$ or $Y_{\tau}^1=Y_{\tau}^2$\,\,$P-$a.s. Therefore $Y^2_{\tau} \geq Y^1_{\tau},$
 for every $\tau\in [S_i, S_{i+1}]$. Then $Y^1\leq Y^2$. The proof of assertion \it{1.} is
 finished.
\\ \it{2.} Let us first point out that if
$$
X_t = X_0 + V_t+\dint_0^t \alpha_s dB_s + \displaystyle{\sum_{0<
s\leq t}} \Delta_s X,
$$
with $X_0\in\R,\,\, V\in {\cal K}^c -{\cal K}^c$ and $\alpha\in{\cal
L}^{2,d}$. Then by using It\^{o}-Tanaka formula, there exists
$\widetilde{l}\in {\cal K}^c$ such that
$$
X_t^+ = X_0^+ +\dint_0^t 1_{\{X_s >0\}}dV_s+\dint_0^t 1_{\{X_s
>0\}}\alpha_s dB_s + \displaystyle{\sum_{0< s\leq t}} \Delta_sX^+
+\widetilde{l}_t.
$$
Henceforth, if $X_t\geq 0$, then
$$
1_{\{X_s=0\}}\alpha_s =0,\,\,dsP(d\omega)-a.e.\,\,\mbox{and}\,\,\,
1_{\{X_{s} =0\}}dV_s = d\widetilde{l}_s \quad\mbox{on}\quad
[0,T],\,\,P-a.s.
$$
Now, by using the above remark to $X_t := Y_t^2 -Y^1_t \geq 0$, then
$1_{\{X_s=0\}}(Z_s^2 -Z_s^1)=0,\,\,dsP(d\omega)-a.e.$ Moreover,
there exists $\widetilde{l}\in {\cal K}^c$ such that
$$
\begin{array}{ll}
- d\widetilde{l}_s &= 1_{\{X_{s}=0\}}(f^2(s,Y_{s}^1,Z_{s}^1)
-f^1(s,Y_{s}^{1},Z_{s}^{1}))ds +
1_{\{{X}_{s}=0\}}(g^2(s,Y_{s}^1)dA_s^2 -g^1(s,Y_{s}^{1})dA_s^1) \\
& +
 1_{\{{X}_{s}=0\}}\bigg(
(d{{K}}_s^{2+,c}-d{{K}}_s^{2-,c})
-(d{{K}}_s^{1+,c}-d{{K}}_s^{1-,c})\bigg).
\end{array}
$$
Therefore
$$
 1_{\{Y_s^1
=Y_s^2\}}d{{K}}_s^{2+, c}\leq 1_{\{ Y_s^1 =Y_s^2\}}(d{{K}}_s^{2-,
c}+ d{{K}}_s^{1+, c}).
$$
Since $dK^{2+,c}\perp dK^{2-,c}$, we get
$$
1_{\{Y_{s}^2 = Y_{s}^1\}}d{{K}}_s^{2+, c}\leq d{{K}}_s^{1+, c}.
$$
Then
$$
1_{\{L_{s}^2 = L_{s}^1\}}d{{K}}_s^{2+, c} =1_{\{L_{s}^2 =
L_{s}^1=Y_{s}^2 = Y_{s}^1\}}d{{K}}_s^{2+, c} \leq d{{K}}_s^{1+, c}.
$$
Let us now compare the discontinuous parts of the reflecting
processes. By formula (\ref{equa1}) and assumption
\noindent$(\bf{D.5})$, we have for each $s\in]0,T]$ such that
$\Delta_sK^{2+}>0$ and $L_{s-}^2 = L_{s-}^1$ (then $L_{s-}^1
=Y_{s-}^1= Y_{s-}^2= L_{s-}^2$)
$$
\begin{array}{ll}
 \Delta_sK^{2+} & = L_{s-}^2- [Y_s^2+h^2(s,  Y_{s-}^2,
 Y_s^2 ) +\Delta_s R^2]
 \\ & =  L_{s-}^1- [Y_s^2+h^2(s, Y_{s-}^1,
 Y_s^2 )+\Delta_s R^2]
 \\ & \leq  L_{s-}^1- [Y_s^1+h^2(s,Y_{s-}^1,
 Y_s^1 )+\Delta_s R^1]
 \\ & \leq  L_{s-}^1- [Y_s^1+h^1(s,Y_{s-}^1,
 Y_s^1 )+\Delta_s R^1]
 \\ & =\Delta_sK^{1+}.
\end{array}
$$
Hence
$$
1_{\{L_{s-}^2 = L_{s-}^1\}}\Delta_sK^{2+}\leq \Delta_sK^{1+}.
$$
Consequently
$$ 1_{\{L^1_{t-} =
L^2_{t-}\}}dK^{2+}_t \leq dK^{1+}_t.
$$
Similarly we have also that
 $$
1_{\{U^1_{t-} = U^2_{t-}\}}dK^{1-}_t \leq dK^{2-}_t.
$$
 This
completes the proof of Theorem \ref{th111}. \eop

%
%
%

\end{document}